\documentclass[a4paper,12pt]{article}
\usepackage{amsmath,amssymb,amsthm}
\usepackage{amscd}
\usepackage[all]{xy}
\newtheorem{dfn}{Definition}[section]
\newtheorem{prop}[dfn]{Proposition}
\newtheorem{thm}[dfn]{Theorem}
\newtheorem{lem}[dfn]{Lemma}
\newtheorem{cor}[dfn]{Corollary}

\begin{document}
\title{The mapping class group and the Meyer function for plane curves}
\author{Yusuke Kuno}
\date{}
\maketitle

\begin{abstract}
For each $d\ge 2$, \textit{the mapping class group for plane curves of degree $d$}
will be defined and it is proved that there exists uniquely \textit{the Meyer function} on this group.
In the case of $d=4$, using our Meyer function,
we can define the local signature for 4-dimensional fiber spaces
whose general fibers are non-hyperelliptic compact
Riemann surfaces of genus 3. Some computations of our local signature will be given.
\end{abstract}

\noindent \textbf{Introduction.}\\
Let $\Sigma_g$ be a closed oriented
$C^{\infty}$-surface of genus $g \ge 0$ and
let $\Gamma_g$ be the mapping class group of $\Sigma_g$, namely the group of
all isotopy classes of orientation preserving diffeomorphisms of $\Sigma_g$.

In \cite{Me} W.\ Meyer discovered and studied a cocycle
$\tau_g \colon \Gamma_g \times \Gamma_g \rightarrow \mathbb{Z}$.
For the sake of the reader a brief definition of $\tau_g$ will be given in Appendix.
This cocycle is called \textit{Meyer's signature cocycle}.
In his paper W.\ Meyer showed that the cohomology class $[\tau_g]\in H^2(\Gamma_g;\mathbb{Z})$ is torsion for
$g=1,2$ and has infinite order for $g\ge 3$, and gave an explicit formula for
the unique $\mathbb{Q}$-valued 1-cochain of $\Gamma_1$ cobounding $\tau_1$
using the Rademacher function (\cite{Me} p.259 Satz 4).
Since the hyperelliptic mapping class group $\Gamma_g^{H}$, a subgroup of $\Gamma_g$, was
shown to be $\mathbb{Q}$-acyclic by F.\ Cohen\cite{FC} and N.\ Kawazumi\cite{K} independently,
it was known to specialists that there
exists the unique 1-cochain of $\Gamma_g^{H}$ cobounding $\tau_g$ restricted to $\Gamma_g^{H}$.
In \cite{E} H.\ Endo directly showed the existence and the uniqueness of such a 1-cochain
$\phi_g^H \colon \Gamma_g^{H} \rightarrow \frac{1}{2g+1}\mathbb{Z}$ using a finite presentation of
$\Gamma_g^{H}$ by J.\ Birman-H.\ Hilden\cite{BH}. He also defined 
the local signature for hyperelliptic fibrations using $\phi_g^H$, and studied the
geometry of hyperelliptic fibrations; for example, he derived a signature formula for
such fibrations over a closed surface. His formula originates from
Y.\ Matsumoto\cite[Theorem 3.3]{Ma} where genus 2 fibrations are discussed.
For the study of the function $\phi_g^H$, see also T.\ Morifuji's paper\cite{Mo}.

The purpose of the present paper is to give another interesting example of these
phenomena; \textit{the Meyer function on the mapping class group for plane curves}.

For $d\ge 2$ a group $\Pi(d)$ and a homomorphism
$\rho \colon \Pi(d)\rightarrow \Gamma_g$,where $g=\frac{1}{2}(d-1)(d-2)$, will be constructed.
The group $\Pi(d)$ can be considered as the fundamental group of the classifying
space for isotopy classes of continuous families of non-singular plane curves of degree $d$;
the precise meaning of this statement will be given in Theorem \ref{thm:6-1} later.

The main results of this paper are Theorem \ref{thm:4-1} and Theorem \ref{thm:4-2}.
As a consequence of them it follows that the pull back $\rho^*[\tau_g]$
vanishes in the rational cohomology $H^2(\Pi(d);\mathbb{Q})$ and there exists the unique 1-cochain
$\phi^d \colon \Pi(d) \rightarrow \mathbb{Q}$ such that $\delta \phi^d=\rho^*\tau_g$.
$\phi^d$ will be called \textit{the Meyer function for plane curves of degree $d$}.

This is similar to the case of $\Gamma_1,\Gamma_2$, and $\Gamma_g^{H}$, but
we remark that the homomorphism $\rho$ seems no more injective nor surjective.
In fact, for $d=4$ we will see in Proposition \ref{prop:6-3} that $\rho$ is
surjective but has non-trivial kernels.
In this sense our result is different from the works of W.\ Meyer and
H.\ Endo where \textit{subgroups} of $\Gamma_g$ are considered.

While they did explicit computations of $\tau_g$ for
certain relators of the mapping class groups to prove the vanishing of $[\tau_g]$,
our method depends on the vanishing of $[\tau_g]$ pulled back to
the cohomology of a fundamental group
of the complement of a hypersurface in a complex vector space, 
which will be stated in Proposition \ref{prop:3-1} and proved
using the definition of Meyer's signature cocycle and the standard argument in
differential topology; the way from Proposition \ref{prop:3-1} to the vanishing of
$[\tau_g]$ pulled back to $H^2(\Pi(d);\mathbb{Q})$ are elementary.
Since this needs no explicit computations of $\tau_g$, we believe that our
method has its own meaning to grasp the conceptual reason of the vanishing of $[\tau_g]$
and can be applied to other cases in the future.

Our study of the vanishing of $\rho^*[\tau_g]\in H^2(\Pi(d);\mathbb{Q})$ has a connection with
localization of the signature of 4-dimensional fiber spaces, that is a recent hot
topic studied in various fields such as topology, algebraic geometry,
and complex analysis (see T.\ Ashikaga-H.\ Endo\cite{AE} and T.\ Ashikaga-K.\ Konno\cite{AK}).

As an application of our study, especially $d=4$, we define the local signature
for the set of all fiber germs of 4-dimensional fiber spaces
whose general fibers are non-hyperelliptic compact
Riemann surfaces of genus 3 by using our 1-cochain $\phi^4$ of $\Pi(4)$.
The fact that any non-hyperelliptic compact Riemann surface of genus 3 can
be realized as a smooth quartic curve in $\mathbb{P}^2$ by the canonical embedding,
is crucial.

In this case of non-hyperelliptic family of genus 3, T.\ Ashikaga-K.\ Konno \cite{AK} and
K.\ Yoshikawa \cite{Y} have already defined local signature independently.
The definition of \cite{AK} is algebro geometric and that of \cite{Y} is complex analytic.
We compute some examples of values of our local signature, defined by topological way,
and observe that they coincide with those computed in \cite{AK} and \cite{Y}.

\section{Definitions}
Throughout this paper, $d$ denotes a fixed integer $\ge 2$.
Let $V^d$ be the complex vector space of homogeneous polynomials of degree $d$
in the determinates $x$,$y$, and $z$, and let $\mathbb{P}(d)=\mathbb{P}(V^d)$ be the
projectivization of $V^d$. By taking the set of monomials
$\{ x^{\ell(k)}y^{m(k)}z^{n(k)} \}_{k=0}^{N}$ of degree $d$, where
$N=\frac{1}{2}(d+2)(d+1)-1$, each element of $V^d$ can be uniquely 
written as the form
$$\Phi=\sum_{k=0}^{N}a_kx^{\ell(k)}y^{m(k)}z^{n(k)},$$
where $a_k\in \mathbb{C}$.
We denote the corresponding homogeneous coordinates of $\mathbb{P}(d)$ by
$[a_0 \colon a_1 \colon \cdots \colon a_N]$. Each element $a\in \mathbb{P}(d)$ determines an algebraic curve
$C_a\subset \mathbb{P}^2$ of degree $d$. Later we also denote by $C_F$ the algebraic curve
defined by $F\in V^d \setminus \{0\}$. We believe this use of notation does not confuse
the reader. Let $D$ be the set of points $a\in \mathbb{P}(d)$
such that the corresponding curve $C_a$ is singular. $D$ is called \textit{the discriminant locus}
and is well-known to be irreducible and of codimension 1.
For a proof, see also the remark after the proof of Proposition \ref{prop:2-1} in this paper.

There is an action of $GL(3;\mathbb{C})$ on $V^d$ given by
$$(A\cdot F)(x,y,z)=F((x,y,z)\cdot {\hskip -3pt}\ ^t{\hskip -1pt}A^{-1}),$$
where $A\in GL(3;\mathbb{C})$ and $F\in V^d.$ Here ${\hskip -3pt}\ ^t{\hskip -1pt}A$
is the transpose of the matrix $A$.  This action induces the action of
$PGL(3)$ on $\mathbb{P}(d)$, $D$, and $\mathbb{P}(d)\setminus D$.

Let $EPGL(3)\rightarrow BPGL(3)$ be the universal principal $PGL(3)$ bundle.
We denote by $\Pi(d)$ the fundamental group of the Borel construction
$(\mathbb{P}(d)\setminus D)_{PGL(3)}=EPGL(3)\times_{PGL(3)}(\mathbb{P}(d)\setminus D)$
and call this group \textit{the mapping class group for plane curves of degree $d$}.

For $(e,a)\in EPGL(3)\times (\mathbb{P}(d)\setminus D)$, we denote by $[e,a]$ the element
of $(\mathbb{P}(d)\setminus D)_{PGL(3)}$ represented by $(e,a)$.
This notation concerning Borel construction will be used several times.

Let $\bar{\mathcal{F}}$ (resp. $\mathcal{F}$) be the hypersurface in
$\mathbb{P}(d) \times \mathbb{P}^2$ (resp. $(\mathbb{P}(d)\setminus D) \times \mathbb{P}^2$)
defined as the zero set of $\Phi$ considered as a bi-homogeneous polynomial in $a_0,\ldots,a_N$
and $x,y,z$.
Then the restriction of the first projection
$p\colon \mathcal{F}\rightarrow \mathbb{P}(d) \setminus D$ is a family of 
non-singular plane curves of degree $d$
whose fiber over $a\in \mathbb{P}(d)\setminus D$ is $C_a$. Since the diagonal action of $PGL(3)$
on $\mathbb{P}(d)\times \mathbb{P}^2$ preserves $\mathcal{F}$ and $p$ is $PGL(3)$-equivariant,
we have a family of Riemann surfaces $p_u\colon \mathcal{F}_{PGL(3)}\rightarrow (\mathbb{P}(d)\setminus D)_{PGL(3)}$.
We denote the topological monodromy (see Appendix) of this family by
$\rho \colon \Pi(d) \rightarrow \Gamma_g$, where $g=\frac{1}{2}(d-1)(d-2)$.
Note that the genus of a non-singular plane curve of degree $d$ is
given by $\frac{1}{2}(d-1)(d-2)$.

In section 4 we will prove that the rational cohomology class
$\rho^*[\tau_g]\in H^2(\Pi(d);\mathbb{Q})$ vanishes and compute the abelianization of $\Pi(d)$.
In section 6 we will prove that the space $(\mathbb{P}(d)\setminus D)_{PGL(3)}$ is the classifying
space of the set of all isotopy classes of continuous families of non-singular plane curves of degree $d$.

\section{The discriminant locus}
In this section we investigate the discriminant locus $D$,
which also can be described in terms of dual variety as follows.
For generality of dual variety, see \cite{GKZ} or \cite{L}. Let $\mathbb{P}(d)^{\vee}$
be the dual projective space of $\mathbb{P}(d)$, i.e., the space of all hyperplanes of $\mathbb{P}(d)$.
We denote by $[\alpha^0 \colon \alpha^1 \colon \cdots \colon \alpha^N]$ the homogeneous coordinates of
$\mathbb{P}(d)^{\vee}$ corresponding to the homogeneous coordinates $[a_0 \colon a_1 \colon \cdots \colon a_N]$
of $\mathbb{P}(d)$; $\alpha=[\alpha^0 \colon \alpha^1 \colon \cdots \colon \alpha^N]$ is the
hypersurface of $\mathbb{P}(d)$ defined by
$$\alpha^0a_0+\alpha^1a_1+\cdots +\alpha^Na_N=0.$$
The Veronese embedding $v\colon \mathbb{P}^2 \rightarrow \mathbb{P}(d)^{\vee}$ is defined by
$$v([x \colon y \colon z])=[x^{\ell(0)}y^{m(0)}z^{n(0)} \colon \cdots \colon x^{\ell(N)}y^{m(N)}z^{n(N)}].$$
Since the dual of $\mathbb{P}(d)^{\vee}$ is canonically isomorphic to $\mathbb{P}(d)$, each element
$a\in \mathbb{P}(d)$ determines the hypersurface of $\mathbb{P}(d)^{\vee}$ which we denote by $H_a$.
We set
$$\mathcal{X}^{\prime}:=\left\{ (a,\alpha)\in \mathbb{P}(d)\times \mathbb{P}(d)^{\vee}
\ ;\ \alpha \in v(\mathbb{P}^2)\  {\rm and}\ H_a \ {\rm is \ tangent \ to}\
v(\mathbb{P}^2)\ {\rm at} \ \alpha \right\}.$$
Then the image of $\mathcal{X}^{\prime}$ by the first projection is just $D$, i.e., $D$ is the dual
variety of $v(\mathbb{P}^2)$.

Let $\mathcal{X}$ be the analytic subset of $\mathbb{P}(d)\times \mathbb{P}^2$ defined by
the equations
$$\Phi=\Phi_x=\Phi_y=\Phi_z=0,$$
where $\Phi_x$ is the partial derivative of $\Phi$ with respect to $x$, etc.
Thus if $(a,p)$ is a point of $\mathcal{X}$, then $a$ is a point of $D$ and $p$ is a
singular point of $C_a$.
Then we see that $\mathcal{X}\rightarrow \mathcal{X}^{\prime},(a,p)\mapsto (a,v(p))$ is an isomorphism.
$\mathcal{X}^{\prime}$ has the structure of fiber bundle over $v(\mathbb{P}^2)$ whose fiber over
$\alpha \in v(\mathbb{P}^2)$ is the set of all hyperplanes in $\mathcal{X}^{\prime}$ tangent to
$v(\mathbb{P}^2)$ at $\alpha$, which is isomorphic to a $(N-3)$-dimensional projective space. 
From this point of view it is clear that $\mathcal{X}$ is non-singular (see also \cite{GKZ}p.30),
but for later consideration we give here an alternative proof using coordinate description.

\begin{prop}
\label{prop:2-1}
$\mathcal{X}$ is non-singular.
\end{prop}

\begin{proof}
Let $(a^0,[x_0 \colon y_0 \colon z_0])$ be a point of $\mathcal{X}$.
We will show $\mathcal{X}$ is non-singular at this point.
Since the action of $PGL(3)$ on $\mathbb{P}(d)\times \mathbb{P}^2$ preserves $\mathcal{X}$,
we may assume that $[x_0 \colon y_0 \colon z_0]=[0 \colon 0 \colon 1]$.
Take a polynomial representative $F\in V^d$ of $a^0$, 
then the coefficient of $z^d$ of $F$ is zero because $[0 \colon 0 \colon 1]\in C_{a^0}$.
Moreover, $F$ cannot be written as the form 
$$F=(\alpha x+\beta y)z^{d-1},$$
where $(\alpha,\beta)\neq (0,0)$ because $[0 \colon 0 \colon 1]$ is a singular point of $C_{a^0}$.
Therefore there is a monomial $x^{\ell(k)}y^{m(k)}z^{n(k)}$ which is different from
$z^d$, $xz^{d-1}$, and $yz^{d-1}$ such that the coefficient of $x^{\ell(k)}y^{m(k)}z^{n(k)}$
of $F$ is not zero. By a rearrangement of indices we may assume that $k=0$ and $a_1$, $a_2$
, and $a_3$ correspond to monomials $z^{d}$, $xz^{d-1}$, and $yz^{d-1}$, respectively.
Then setting $a_0=1$ and $z=1$, we have an inhomogeneous coordinates
$(a_1,\ldots,a_N,x,y)$ of $\mathbb{P}(d)\times \mathbb{P}^2$ near $(a^0,[0 \colon 0 \colon 1])$.
In this local coordinate system $\mathcal{X}$ is defined by the equations
$$\Psi=\Psi_x=\Psi_y=0,$$
where $\Psi=\Phi(1,a_1,\ldots,a_N,x,y,1)$.

Now the Jacobian matrix of $(\Psi,\Psi_x,\Psi_y)$ at $(a^0,[0 \colon 0 \colon 1])$ is 
\begin{eqnarray*}
J &=& \left( \begin{array}{ccccccc}
\Psi_{a_1} & \Psi_{a_2} & \Psi_{a_3} & \cdots & \Psi_x & \Psi_y \\
\Psi_{x,a_1} & \Psi_{x,a_2} & \Psi_{x,a_3} & \cdots & \Psi_{xx} & \Psi_{xy} \\
\Psi_{y,a_1} & \Psi_{y,a_2} & \Psi_{y,a_3} & \cdots & \Psi_{yx} & \Psi_{yy} \\
\end{array} \right) \\
 &=& \left( \begin{array}{ccccccc}
1 & 0 & 0 & \cdots & 0 & 0 & 0 \\
0 & 1 & 0 & \cdots & 0 & \Psi_{xx} & \Psi_{xy} \\
0 & 0 & 1 & \cdots & 0 & \Psi_{yx} & \Psi_{yy} \\
\end{array} \right),
\end{eqnarray*}
we see that the rank of $J$ is 3. This shows that $\mathcal{X}$ is non-singular at
$(a^0,[0 \colon 0 \colon 1])$.
\end{proof}

Let $\pi \colon \mathcal{X} \rightarrow D \subset \mathbb{P}(d)$ be the first projection.
The above proof shows that $(a^0,[0\colon 0\colon 1])$ is a regular point of $\pi$ if and
only if
$${\rm det}\left( \begin{array}{cc}
\Psi_{xx} & \Psi_{xy} \\
\Psi_{yx} & \Psi_{yy} \\
\end{array} \right)\neq 0$$
at $(a^0,[0\colon 0\colon 1])$.
By an argument like the Morse lemma, we can take a coordinate system $(X,Y)$ of $\mathbb{P}^2$
centered at $[0 \colon 0\colon 1]$ such that $C_{a^0}$ is locally given by the equation $X^2+Y^2=0$.
Thus $[0\colon 0\colon 1]$ is a nodal singularity. This holds for other points of $\mathcal{X}$;
$(a,p)\in \mathcal{X}$ is a regular point of $\pi$ if and only if $p$ is a nodal singularity of $C_a$.

Let $E$ be the union of singular points of $D$ and the $\pi$-image of critical points of $\pi$. 
$E$ is a proper analytic subset of $D$ by Sard's theorem.

Here we give a short proof that $D$ is irreducible and of codimension 1.
At first, $\mathcal{X}\cong \mathcal{X}^{\prime}$ is non-singular and connected hence irreducible. 
Therefore, $D=\pi(\mathcal{X})$ is also irreducible.
On the other hand $D$ is at most $N-1$ dimensional because $D$ is a proper analytic subset of $\mathbb{P}(d)$.
Let $a$ be a point of $D\setminus E$ and take a point $(a,p)$ in the fiber $\pi^{-1}(a)$.
Then $D$ is smooth around $a$ and the differential of $\pi$ at $(a,p)$ is of maximal rank $N-1$.
This shows $D$ is indeed $N-1$ dimensional. Note that $E$ is at most $N-2$ dimensional.

In the next lemma we shall describe the hyperplane of $\mathbb{P}(d)$ tangent to
$D$ at a point in $D\setminus E$.

\begin{lem}
\label{lem:2-2}
Let $(a^0,[x_0 \colon y_0 \colon z_0])$ be a point of $\mathcal{X}$ and suppose that
$a^0 \in D\setminus E$. Then the hyperplane $T_{a^0}$ tangent to $D$ at $a^0$ is given by
$$T_{a^0}=\left\{ [\xi_0 \colon \xi_1 \colon \cdots \colon \xi_N]\in \mathbb{P}(d) \ ;
\ \sum_{k=0}^{N}\xi_kx_0^{\ell(k)}y_0^{m(k)}z_0^{n(k)}=0 \right\}.$$
Moreover, $[x_0 \colon y_0 \colon z_0]$ is the unique singular point of $C_{a^0}$.
\end{lem}
\begin{proof}
To prove the first part,
we may assume $a_0^0=z_0=1$ and take an inhomogeneous coordinate system
$(a_1,\ldots,a_N,x,y)$ of $\mathbb{P}(d)\times \mathbb{P}^2$ near $(a^0,[x_0 \colon y_0 \colon 1])$.
Since $a^0$ is a non-singular point of $D$ and $(a^0,[x_0 \colon y_0 \colon 1])$ is a regular
point of $\pi$, we have
$T_{a^0}D=\tilde{\pi}_*(T_{(a^0,[x_0 \colon y_0 \colon 1])}\mathcal{X}),$
where $\tilde{\pi}_*  \colon T_{(a^0,[x_0 \colon y_0 \colon 1])}(\mathbb{P}(d)\times \mathbb{P}^2) \rightarrow
T_{a^0}\mathbb{P}(d)$ is the differential of the first projection
$\tilde{\pi} \colon \mathbb{P}(d)\times \mathbb{P}^2 \rightarrow \mathbb{P}(d)$ and we regard
$T_{(a^0,[x_0 \colon y_0 \colon 1])}\mathcal{X}$ (resp. $T_{a^0}D$) as the subspace of
$T_{(a^0,[x_0 \colon y_0 \colon 1])}(\mathbb{P}(d)\times \mathbb{P}^2)$ (resp. $T_{a^0}\mathbb{P}(d)$).

Now the Jacobian matrix $J$ appeared in the proof of Proposition \ref{prop:2-1}
has the form
$$J = \left( \begin{array}{ccccc}
x_0^{\ell(1)}y_0^{m(1)} & \cdots & x_0^{\ell(N)}y_0^{m(N)}& 0 & 0 \\
* & \cdots & * & \Psi_{xx} & \Psi_{xy} \\
* & \cdots & * & \Psi_{yx} & \Psi_{yy} \\
\end{array} \right)$$
at $(a^0,[x_0\colon y_0\colon 1])$. The rank of this matrix is 3, because
${\rm det}\left( \begin{array}{cc}
\Psi_{xx} & \Psi_{xy} \\
\Psi_{yx} & \Psi_{yy} \\
\end{array} \right) \neq 0$ at $(a^0,[x_0\colon y_0\colon 1])$
by $a^0 \notin E$ and there is an index $i$ such that $x_0^{\ell(i)}y_0^{m(i)}\neq 0$.

Therefore
$$T_{(a^0,[x_0 \colon y_0 \colon 1])}\mathcal{X}=\left\{ \sum_{k=1}^{N}\xi_k\frac{\partial}{\partial a_k}
+\xi_{N+1}\frac{\partial}{\partial x}+\xi_{N+2}\frac{\partial}{\partial y} \ ;
\ J\left(\begin{array}{c}
\xi_1 \\
\vdots \\
\xi_{N+2} \\
\end{array} \right)=0 \right\}$$
and
$$T_{a^0}D=\tilde{\pi}_*(T_{(a^0,[x_0 \colon y_0 \colon 1])}\mathcal{X})=
\left\{ \sum_{k=1}^{N}\xi_k\frac{\partial}{\partial a_k}
\ ;\ \sum_{k=1}^N\xi_kx_0^{\ell(k)}y_0^{m(k)}=0 \right\}.$$
Interpretting this equation in terms of homogeneous coordinates of $\mathbb{P}(d)$,
we obtain the desired description of $T_{a^0}$.
The latter statement of the lemma follows from the form of $T_{a^0}$ just proved and
the injectivity of the Veronese embedding.
\end{proof}

Combining the remark after the proof of Proposition \ref{prop:2-1},
we can say more about the curve $C_{a^0}$:

\begin{lem}
\label{lem:2-3}
Let $a^0 \in D \setminus E$ and $[x_0 \colon y_0 \colon z_0]$ be as in Lemma \ref{lem:2-2}.
Then $[x_0 \colon y_0 \colon z_0]$ is a nodal singularity of $C_{a^0}$, and $C_{a^0}$ is
irreducible except for $d=2$.
Thus if $d\ge 3$ the topological type of $C_{a^0}$ is Lefschetz singular fiber of type I,
that is obtained by pinching a non-separating
simple closed curve on $\Sigma_g$ into a point.
\end{lem}
\begin{proof}We only have to show the irreducibility of $C_{a^0}$
for $d\ge 3$. If $C_{a^0}$ is reducible it has two irreducible components $C_1$ and $C_2$
with degrees $d_1$ and $d_2$, and they intersect transversely at one point.
We have $d_1d_2=1$ by B\'ezout's theorem, but this contradicts to $d_1+d_2=d\ge 3$.
\end{proof}

The projective space $\mathbb{P}(d)$ can be regarded as the set of all complex lines through
the origin in $V^d$. Let $\widetilde{D}$ (resp. $\widetilde{E}$) be the union of all lines in $D$ (resp. $E$).
In the coordinate system $(a_0,\ldots,a_N)$ of $V^d$, the tangent space of $\widetilde{D}$ at
$F\in \widetilde{D}\setminus \widetilde{E}$ is given by
$$T_F\widetilde{D}=\left\{ \sum_{k=0}^N\xi_k\frac{\partial}{\partial a_k}
\ ;\ \sum_{k=0}^N\xi_kx_0^{\ell(k)}y_0^{m(k)}z_0^{n(k)}=0 \right\},$$
where $[x_0 \colon y_0 \colon z_0]$ is the singular point of $C_F$.
This follows from Lemma \ref{lem:2-2}.

We shall prove a useful lemma which will be used in the next two sections.
Let $\widetilde{\mathcal{F}}$ be the family of algebraic curves over
$V^d \setminus\{0\}$ defined as in the case of $\bar{\mathcal{F}}$ over $\mathbb{P}(d)$.

\begin{lem}
\label{lem:2-4}
Let $B$ be a $C^{\infty}$-manifold of dimension $s\ge 2$ and $j \colon B\rightarrow V^d$ a
$C^{\infty}$-map such that
$j(B)\subset V^d \setminus \widetilde{E}$ and $j$ is transverse to $\widetilde{D}$.
Then the total space $j^*\widetilde{\mathcal{F}}$ of the pull back of
the family $\widetilde{\mathcal{F}}$ by $j$ is a $C^{\infty}$-manifold.
\end{lem}
\begin{proof}
$j^*\widetilde{\mathcal{F}}$ is given by 
$$j^*\widetilde{\mathcal{F}}=\left\{ (b,p)\in B\times \mathbb{P}^2\ ;\ \Phi(j(b),p)=0 \right\}$$
and it is easy to see that if $(b^0,p_0)\in j^*\widetilde{\mathcal{F}}$ and 
$p_0$ is a smooth point of $C_{j(b^0)}$ then $j^*\widetilde{\mathcal{F}}$ is smooth at $(b^0,p_0)$.

Suppose $(b^0,p_0)\in j^*\widetilde{\mathcal{F}}$ and $p_0=[x_0 \colon y_0 \colon z_0]$
is the singular point of $C_{j(b^0)}$.
Note that we have $j(b^0)\in \widetilde{D} \setminus \widetilde{E}$.
Let $(j_0,j_1,\ldots,j_N)$ denote the $N+1$-tuples of smooth functions on $B$ determined by
$j$ and the coordinate system $(a_0,a_1,\ldots,a_N)$ of $V^d$. 
By the assumption of transversality and the description of $T_{j(b^0)}\widetilde{D}$ given above, 
we can choose a suitable local coordinate system $(b_1,\ldots,b_s)$ of $B$ around $b_0$ such that complex numbers
$$\sum_{k=0}^N\frac{\partial j_k}{\partial b_1}(b^0)x_0^{\ell(k)}y_0^{m(k)}z_0^{n(k)} \ {\rm and}\ 
\sum_{k=0}^N\frac{\partial j_k}{\partial b_2}(b^0)x_0^{\ell(k)}y_0^{m(k)}z_0^{n(k)}$$
are linearly independent over the real numbers. From this we can conclude that
$j^*\widetilde{\mathcal{F}}$ is smooth at $(b^0,p_0)$. This completes the proof.
\end{proof}

We remark that in holomorphic category one can say more;
if $B$ is a complex manifold of complex dimension $\ge 1$ and $j$ is holomorphic,
$j^*\widetilde{\mathcal{F}}$ has a complex
structure as a hypersurface in $B \times \mathbb{P}^2$.

\section{The 1-cochain of $\pi_1(V^d \setminus \widetilde{D})$}
Let $\chi_1 \colon \pi_1(V^d \setminus \widetilde{D}) \rightarrow \pi_1(\mathbb{P}(d)\setminus D)$
be the homomorphism induced by the projection map
$V^d \setminus \widetilde{D} \rightarrow \mathbb{P}(d)\setminus D$
and let $\chi_2 \colon \pi_1(\mathbb{P}(d)\setminus D) \rightarrow \Pi(d)$ be the homomorphism
induced by the inclusion map $\mathbb{P}(d)\setminus D \rightarrow (\mathbb{P}(d)\setminus D)_{PGL(3)},
a \mapsto [e_0,a]$ where $e_0$ is the base point of $EPGL(3)$.
We set $\chi:=\chi_2 \circ \chi_1$ and $\tilde{\rho}:=\rho \circ \chi$.
Then $\tilde{\rho} \colon \pi_1(V^d \setminus \widetilde{D})\rightarrow \Gamma_g$ is the topological
monodromy of the family over $V^d \setminus \widetilde{D}$
defined as in the case of $\mathcal{F}\rightarrow \mathbb{P}(d)\setminus D$.

In this section, we shall construct a 1-cochain
$c \colon \pi_1(V^d \setminus \widetilde{D})\rightarrow \mathbb{Z}$ and prove that $\delta c=\tilde{\rho}^*\tau_g$.
The key is that $V^d \setminus \widetilde{E}$ is 2-connected, which follows from the
fact that the complex codimension of $\widetilde{E}$ in $V^d$ is $\ge 2$.
All of the spaces that we consider in this section as well as all of the maps are based.

We regard the circle $S^1$ as the boundary of the unit disk $D^2$ in $\mathbb{R}^2$.
$D^2$ has the natural orientation induced by that of $\mathbb{R}^2$ and this induces 
the orientation of $S^1$ by counter clockwise manner.
Let $\ell \colon S^1 \rightarrow V^d \setminus \widetilde{D}$
be a $C^{\infty}$-map. Since $V^d \setminus \widetilde{E}$ is simply connected
we can extend $\ell$ to a $C^{\infty}$-map $\tilde{\ell} \colon D^2 \rightarrow V^d \setminus \widetilde{E}$.
We may assume that $\tilde{\ell}$ is transverse to $\widetilde{D}$.
By Lemma \ref{lem:2-4} $\tilde{\ell}^*\widetilde{\mathcal{F}}$ is a compact 4-dimensional
$C^{\infty}$-manifold with boundary
and has the natural orientation induced by the orientation of $D^2$ and that of the fibers,
which have the natural orientations as compact Riemann surfaces.
Set 
$$c([\ell]):={\rm Sign}(\tilde{\ell}^*\widetilde{\mathcal{F}}),$$ where $[\ell]$ denotes
the element of $\pi_1(V^d \setminus \widetilde{D})$ represented by $\ell$ and the right hand side
is the signature of $\tilde{\ell}^*\widetilde{\mathcal{F}}$.

\begin{prop}
\label{prop:3-1}
The above definition of $c$ is well defined and $\delta c=\tilde{\rho}^*\tau_g$, i.e.,
$c$ is a cobounding cochain for $\tilde{\rho}^*\tau_g$.
\end{prop}
\begin{proof}
We first show that $c$ is well defined. Let $\ell_0$ and
$\ell_1$ are $C^{\infty}$-maps from $S^1$ to $V^d \setminus \widetilde{D}$, and suppose
that they represent the same element of $\pi_1(V^d \setminus \widetilde{D})$.
Then there exists a $C^{\infty}$-homotopy $H \colon S^1 \times [0,1] \rightarrow V^d \setminus \widetilde{D}$
such that $H(\cdot,0)=\ell_0$ and $H(\cdot,1)=\ell_1$.
 
Regard the 2-sphere $S^2$ as the annulus $S^1 \times [0,1]$ with two copies of
$D^2$ attached along its two boundary circles $S^1 \times \{0\}$ and $S^1 \times \{1\}$.
We denote by $D_0^2$ one of copies of $D^2$ attached to $S^1 \times \{0\}$ and $D_1^2$ the other.
Using some extensions $\tilde{\ell_i} \colon D_i^2 \rightarrow V^d \setminus \widetilde{E}$
of $\ell_i$ for $i=0$ and $1$,
$H$ extends to a $C^{\infty}$-map $\widetilde{H} \colon S^2 \rightarrow V^d \setminus \widetilde{E}$.
We introduce the orientation of $S^2$ so that the inclusion $D_0^2 \hookrightarrow S^2$ is
orientation preserving. Thus the other inclusion $D_1^2 \hookrightarrow S^2$ is
orientation reversing.

Since  $\pi_2(V^d \setminus \widetilde{E})=0$, we can extend $\widetilde{H}$ to a $C^{\infty}$-map
$\bar{H} \colon D^3 \rightarrow V^d \setminus \widetilde{E}$ transverse to
$\widetilde{D}\setminus \widetilde{E}$.
Then $\bar{H}^*\widetilde{\mathcal{F}}$ is a
$C^{\infty}$-manifold with boundary $\widetilde{H}^*\widetilde{\mathcal{F}}$.
Since the signature of the boundary of a manifold is zero, we have by the Novikov additivity of the signature
$${\rm Sign}(\tilde{\ell_0}^*\widetilde{\mathcal{F}})-{\rm Sign}(\tilde{\ell_1}^*\widetilde{\mathcal{F}})=0,$$
so $c$ is well defined.

We next show the latter part.
Let $\ell_0$ and $\ell_1$ be $C^{\infty}$-maps from $S^1$ to $V^d \setminus \widetilde{D}$.
We will show
\begin{equation}
\label{eq:1}
c([\ell_0])+c([\ell_1])-c([\ell_0][\ell_1])=\tilde{\rho}^*\tau_g([\ell_0],[\ell_1]).
\end{equation}
Let $P$ denote the pair of pants; this is the 2-sphere $S^2$ with the interior of the three
disjoint closed disks removed.
We also choose two of three boundary components of $P$ and denote them by
$S_0^1$ and $S_1^1$, respectively. $S_0^1$ and $S_1^1$ have the natural orientations induced by that 
of $P$ and can be naturally identified with $S^1$. 
Since $P$ has the homotopy type of the bouquet $S^1 \vee S^1$,
we can construct a $C^{\infty}$-map $L \colon P \rightarrow V^d \setminus \widetilde{D}$ such that
the restriction of $L$ to $S_0^1$ (resp. $S_1^1$) are exactly same as $\ell_0$ (resp. $\ell_1$).

We notice that the restriction of $L$ to the remaining boundary component of $P$ with the natural orientation
is homotopic to the inverse of the composition loop $\ell_0 \cdot \ell_1$. We also have
${\rm Sign}(L^*\widetilde{\mathcal{F}})=-\tilde{\rho}^*\tau_g([\ell_0],[\ell_1])$
by the definition of Meyer's signature cocycle $\tau_g$.
Using some extensions $\tilde{\ell_i}$ of $\ell_i$ for $i=0$ and $1$, and an extension 
$\widetilde{\ell_0 \cdot \ell_1}$ of $\ell_0 \cdot \ell_1$,
$L$ extends to a $C^{\infty}$-map $\widetilde{L} \colon S^2 \rightarrow V^d \setminus \widetilde{E}$.
Moreover $\widetilde{L}$ extends to a map $\bar{L} \colon D^3 \rightarrow V^d \setminus \widetilde{E}$
transverse to $\widetilde{D}$. 
We have ${\rm Sign}(\widetilde{L}^*\widetilde{\mathcal{F}})=0$ since $\widetilde{L}^*{\mathcal{F}}$
is the boundary of $\bar{L}^*{\mathcal{F}}$ hence we obtain by the Novikov additivity 
$$0={\rm Sign}(\widetilde{L}^*\widetilde{\mathcal{F}})={\rm Sign}(\tilde{\ell_0}^*\widetilde{\mathcal{F}})
+{\rm Sign}(\tilde{\ell_1}^*\widetilde{\mathcal{F}})
-{\rm Sign}(\widetilde{\ell_0 \cdot \ell_1}^*\widetilde{\mathcal{F}})+{\rm Sign}(L^*\widetilde{\mathcal{F}}),$$
that is just the equation (\ref{eq:1}).
\end{proof}

\section{Main theorems}
In this section we shall state and prove the main results of this paper.
In section 1 we defined the group $\Pi(d)$ and the homomorphism $\rho \colon \Pi(d) \rightarrow \Gamma_g$.
\begin{thm}$\rho^*[\tau_g]=0 \in H^2(\Pi(d);\mathbb{Q})$.
\label{thm:4-1}
\end{thm}

\begin{thm}The first homology group of $\Pi(d)$ is given as follows:
$$H_1(\Pi(d);\mathbb{Z})=
\begin{cases}
\mathbb{Z}/3(d-1)^2\mathbb{Z} & {\rm if}\  d\equiv 0 \ {\rm mod}\ 3, \\
\mathbb{Z}/(d-1)^2\mathbb{Z} & {\rm if}\ d\equiv 1 \ {\rm or}\ 2 \ {\rm mod}\ 3. 
\end{cases}$$
In particular, we have $H^1(\Pi(d);\mathbb{Q})=0$.
\label{thm:4-2}
\end{thm}

As an immediate consequence of these theorems, it follows that there exists the unique 1-cochain 
$\phi^d \colon \Pi(d) \rightarrow \mathbb{Q}$ such that $\delta \phi^d=\rho^*\tau_g$. We will call $\phi^d$
\textit{the Meyer function for plane curves of degree $d$}.

The rest of this section will be devoted to the proof of these theorems.
In Proposition \ref{prop:3-1} we have showed that
$\tilde{\rho}^*[\tau_g]=0 \in H^2(\pi_1(V^d \setminus \widetilde{D});\mathbb{Z})$.
Thus Theorem \ref{thm:4-1} follows from the following:
\begin{lem}
The homomorphism 
$$\chi^*:H^2(\Pi(d);\mathbb{Q}) \rightarrow H^2(\pi_1(V^d \setminus \widetilde{D});\mathbb{Q})$$
induced by $\chi$, introduced in section 3, is injective.
\end{lem}
\begin{proof}
Recall that $\chi$ is the composition of
$\chi_1$ and $\chi_2$. We first consider $\chi_1$.
Let $\xi \in H^2(\mathbb{P}(d);\mathbb{Q})$ denote the first Chern class of the principal $\mathbb{C}^*$ bundle
$V^d \setminus \{0\} \rightarrow \mathbb{P}(d)$. Then the restriction of $\xi$ to $\mathbb{P}(d)\setminus D$
is zero, for the first Chern class $c_1([D]) \in H^2(\mathbb{P}(d);\mathbb{Q})$ of the line bundle $[D]$
determined by the divisor $D$ of $\mathbb{P}(d)$ is a multiple of $\xi$ and of course the restriction
of $c_1([D])$ to $\mathbb{P}(d)\setminus D$ is zero.

By the Gysin sequence
$$H^0(\mathbb{P}(d)\setminus D;\mathbb{Q}) \stackrel{\cup \xi}{\rightarrow}H^2(\mathbb{P}(d)\setminus D;\mathbb{Q})
\rightarrow H^2(V^d \setminus \widetilde{D};\mathbb{Q})$$
of the principal $\mathbb{C}^*$ bundle
$V^d \setminus \widetilde{D} \rightarrow \mathbb{P}(d) \setminus D$
we see that $H^2(\mathbb{P}(d)\setminus D;\mathbb{Q}) \rightarrow H^2(V^d \setminus \widetilde{D};\mathbb{Q})$
is injective.
Therefore
$$\chi_1^* \colon H^2(\pi_1(\mathbb{P}(d)\setminus D);\mathbb{Q}) \rightarrow
H^2(\pi_1(V^d \setminus \widetilde{D});\mathbb{Q})$$
is also injective.

We next consider $\chi_2$. By the homotopy exact sequence of the $\mathbb{P}(d)\setminus D$ bundle
$(\mathbb{P}(d)\setminus D)_{PGL(3)} \rightarrow BPGL(3), [e,a]\mapsto \varpi(e)$
where $\varpi$ denotes the projection map $EPGL(3)\rightarrow BPGL(3)$, we have an exact sequence
\begin{equation}
\mathbb{Z}/3\mathbb{Z} \cong \pi_2(BPGL(3)) \rightarrow \pi_1(\mathbb{P}(d)\setminus D)
\stackrel{\chi_2}{\rightarrow} \Pi(d) \rightarrow 1.
\label{eq:2}
\end{equation}
This implies that
$$\chi_2^* \colon H^2(\Pi(d);\mathbb{Q}) \rightarrow H^2(\pi_1(\mathbb{P}(d)\setminus D);\mathbb{Q})$$
is isomorphic.
Since $\chi^*=\chi_1^* \circ \chi_2^*$, the lemma follows.
\end{proof}

We next proceed to Theorem \ref{thm:4-2}. In the following we consider (co)homology with
coefficients in $\mathbb{Z}$. We need the following two lemmas:

\begin{lem}
\label{lem:4-4}
Let $a^0 \in \mathbb{P}(d)\setminus D$ and denote by $\mathbb{P}$ the set of all
projective lines in $\mathbb{P}(d)$ through $a^0$. Then there exist a non-empty Zariski open
subset $U \subset \mathbb{P}$ such that each element of $U$ does not meet $E$
and is transverse to $D$.
\end{lem}
\begin{proof}
Consider the projection with center $a^0$
$$f\colon D \rightarrow \mathbb{P},\  f(a)={\rm the \ line \ through}\ a^0\ {\rm and}\ a.$$
Note that for $a\in D\setminus E$, $f$ is critical at $a$ if and only if $f(a)$ is contained
in the hyperplane $T_a$ appeared in Lemma \ref{lem:2-2}, namely $f(a)$ is not transverse
to $D$ at $a$.

$\mathbb{P}$ is a $(N-1)$-dimensional projective space and
$f(E)$ is a proper algebraic set in $\mathbb{P}$ since ${\rm dim}E \le N-2$.
Let $K$ denote the set of all critical values of $f\circ \pi \colon \mathcal{X}\rightarrow \mathbb{P}$.
$K$ contains all critical values of $f|_{D\setminus E}$ since
$\pi|_{\pi^{-1}(D\setminus E)}\colon \pi^{-1}(D\setminus E) \rightarrow D\setminus E$ is
biholomorphic by Lemma \ref{lem:2-2}, and
is algebraic and proper because $K$ is nowhere dense in $\mathbb{P}$ by Sard's theorem.
Therefore if we set
$$U:= \mathbb{P}\setminus (f(E)\cup K),$$
$U$ has the desired property.
\end{proof}

\begin{lem}
\label{lem:4-5}
Let $a^0$ and $U$ be as in Lemma \ref{lem:4-4}. For each $Q \in U$ the invariants of
the complex surface $M=\{ (a,p)\in Q \times \mathbb{P}^2 \ ;\ p \in C_{a} \}$
is given as follows:
$$c_1^2(M)=-d^2+9,\ c_2(M)=d^2+3,\  {\rm Sign}(M)=1-d^2.$$
\end{lem}
\begin{proof}
Since $Q\cong \mathbb{P}^1$ we can regard $M$ as a smooth
hypersurface in $\mathbb{P}^1\times \mathbb{P}^2$ determined by a $(1,d)$ homogeneous polynomial.
For $i=1$ and $i=2$ respectively, we denote by $\xi_i\in H^2(\mathbb{P}^1 \times \mathbb{P}^2;\mathbb{Z})$
the pull back of the first Chern class of $\mathcal{O}(1)$ by $H^2(\mathbb{P}^i;\mathbb{Z}) \rightarrow
H^2(\mathbb{P}^1 \times \mathbb{P}^2;\mathbb{Z})$ induced by the projection $\mathbb{P}^1 \times \mathbb{P}^2
\rightarrow \mathbb{P}^i$. Here $\mathcal{O}(1)$ denotes the dual of the tautological line bundle over
$\mathbb{P}^i$.
The first Chern class of the line bundle $[M]$ determined by the divisor $M$ of
$\mathbb{P}^1\times \mathbb{P}^2$ is $c_1([M])=\xi_1+d\xi_2$.
Therefore by the adjunction formula, the first Chern class of $M$ is
\begin{eqnarray*}
c_1(M) &=& (c_1(\mathbb{P}^1\times \mathbb{P}^2)-c_1([M]))|_M \\
&=& \left(2\xi_1+3\xi_2-(\xi_1+d\xi_2)\right)|_M \\
&=& \left(\xi_1+(3-d)\xi_2\right)|_M.
\end{eqnarray*}
Then the Chern number $c_1^2(M)$ is computed as follows:
\begin{eqnarray*}
c_1^2(M) &=& \langle c_1(M)^2,\mu_M \rangle \\
&=& \langle (\xi_1+(3-d)\xi_2)^2c_1([M]),\mu_{\mathbb{P}^1\times \mathbb{P}^2} \rangle \\
&=& \langle (\xi_1+(3-d)\xi_2)^2(\xi_1+d\xi_2),\mu_{\mathbb{P}^1\times \mathbb{P}^2} \rangle \\
&=& \langle (-d^2+9)\xi_1\xi_2^2,\mu_{\mathbb{P}^1\times \mathbb{P}^2} \rangle \\
&=& -d^2+9.
\end{eqnarray*}
Here $\mu_M$ (resp. $\mu_{\mathbb{P}^1\times \mathbb{P}^2}$) denotes the fundamental homology class of
$M$ (resp. $\mathbb{P}^1\times \mathbb{P}^2$) and
$\langle - ,- \rangle$ denotes the Kronecker pairing between cohomology and homology.
We next compute $c_2(M)$. Again by the adjunction formula, the second Chern class of $M$ is
\begin{eqnarray*}
c_2(M) &=& c_2(\mathbb{P}^1\times \mathbb{P}^2)|_M-c_1(M)\cdot c_1([M])|_M \\
&=& \left( 3\xi_2^2+6\xi_1\xi_2-(\xi_1+(3-d)\xi_2)(\xi_1+d\xi_2)\right)|_M \\
&=& \left( 3\xi_1\xi_2+(d^2-3d+3)\xi_2^2 \right)|_M,
\end{eqnarray*}
and the Chern number which will also be denoted by $c_2(M)$ is 
\begin{eqnarray*}
c_2(M) &=& \langle c_2(M),\mu_M \rangle \\
&=& \langle (3\xi_1\xi_2+(d^2-3d+3)\xi_2^2)c_1([M]),\mu_{\mathbb{P}^1\times \mathbb{P}^2} \rangle \\
&=& \langle (3\xi_1\xi_2+(d^2-3d+3)\xi_2^2)(\xi_1+d\xi_2),\mu_{\mathbb{P}^1\times \mathbb{P}^2} \rangle \\
&=& \langle(d^2+3)\xi_1\xi_2^2,\mu_{\mathbb{P}^1\times \mathbb{P}^2} \rangle \\
&=& d^2+3.
\end{eqnarray*}
Finally by the Hirzebruch signature theorem we have ${\rm Sign}(M)=\frac{1}{3}(c_1^2(M)-2c_2(M))=1-d^2$.
\end{proof}
 
Let $a^0$ and $Q \in U$ be as in Lemma \ref{lem:4-5}. Using the above two lemmas 
we can compute the first homology group of $\pi_1(\mathbb{P}(d) \setminus D)$:
\begin{prop}
\label{prop:4-6}
$H_1(\pi_1(\mathbb{P}(d) \setminus D);\mathbb{Z})=\mathbb{Z}/3(d-1)^2\mathbb{Z}.$
\end{prop}
\begin{proof}
The first projection $g\colon M \rightarrow Q$ is
a family of algebraic curves, whose all singular fibers are of type I by Lemma \ref{lem:2-3}.
Since the Euler contribution (see \cite{BPV}p.118, (11.4)Proposition) of a singular fiber of type I is +1,
the number of singular fibers of $g\colon M \rightarrow Q \cong \mathbb{P}^1$ is 
\begin{equation}
\label{eq:3}
c_2(M)-2(2-2g)=d^2+3-2\left(2-2\cdot \frac{1}{2}(d-1)(d-2)\right)=3(d-1)^2.
\end{equation} 
Now consider the following commutative diagram:
$$\xymatrix{
\mathbb{Z} \cong H_2(\mathbb{P}(d)) \ar[d]_{\cong} \ar[r] &
H_2(\mathbb{P}(d),\mathbb{P}(d) \setminus D) \ar[d]^{\cong} \ar[r]
& H_1(\mathbb{P}(d) \setminus D) \ar[r] & 0 \\
H^{2N-2}(\mathbb{P}(d)) \ar[r]^{\iota^*} & H^{2N-2}(D)\cong \mathbb{Z}. \\
}
$$
Here the vertical isomorphisms are Poincar\'e duality and the first horizontal
sequence is a part of the homology exact sequence of the pair $(\mathbb{P}(d),\mathbb{P}(d) \setminus D)$,
and $\iota^*$ is induced by the inclusion $D \hookrightarrow \mathbb{P}(d)$.
For a generator of $H_2(\mathbb{P}(d))$ we can choose $[Q]$. We can
conclude this generator is mapped to $3(d-1)^2$ times a generator of $H^{2N-2}(D)$ in the above diagram,
because (\ref{eq:3}) shows that $Q$ and $D$ intersect transversally in $3(d-1)^2$ points.
This completes the proof, since $H_1(\mathbb{P}(d) \setminus D)\cong H_1(\pi_1(\mathbb{P}(d) \setminus D);\mathbb{Z})$
is isomorphic to the cokernel of
$H_2(\mathbb{P}(d))\rightarrow H_2(\mathbb{P}(d),\mathbb{P}(d) \setminus D)$.
\end{proof}

Now we start the proof of Theorem \ref{thm:4-2}.
Let $F_0 \in V^d \setminus \widetilde{D}$ be a base point
and $a^0$ the image of $F_0$ under the map $V^d \setminus \widetilde{D} \rightarrow \mathbb{P}(d) \setminus D$.
We consider the maps $\lambda \colon GL(3;\mathbb{C})\rightarrow V^d \setminus \widetilde{D}, A \mapsto A\cdot F_0$ and
$\bar{\lambda} \colon PGL(3)\rightarrow \mathbb{P}(d) \setminus D, \bar{A} \mapsto \bar{A}\cdot a^0$.
Since the isomorphism $\pi_2(BPGL(3))\cong \pi_1(PGL(3))$ induced by the homotopy exact sequence of the
universal $PGL(3)$ bundle $EPGL(3)\rightarrow BPGL(3)$ is compatible with (\ref{eq:2}) and
$\bar{\lambda}_* \colon \pi_1(PGL(3))\rightarrow \pi_1(\mathbb{P}(d) \setminus D)$, we have an exact sequence
of group homology
$$\mathbb{Z}/3\mathbb{Z} \cong H_1(\pi_1(PGL(3))) \stackrel{\bar{\lambda}_*}{\rightarrow}
H_1(\pi_1(\mathbb{P}(d)\setminus D))\cong \mathbb{Z}/3(d-1)^2\mathbb{Z}
\stackrel{{\chi_2}_*}{\rightarrow} H_1(\Pi(d)) \rightarrow 0.$$
Therefore we must compute the map $\bar{\lambda}_*$ to determine $H_1(\Pi(d))$.

For this purpose, we consider the following exact sequence
$$\mathbb{Z} \cong H_1(\mathbb{C}^*)\rightarrow H_1(V^d \setminus \widetilde{D}) \rightarrow H_1(\mathbb{P}(d)\setminus D)
\rightarrow 0$$
induced by a part of the homotopy exact sequence of the principal $\mathbb{C}^*$ bundle
$V^d \setminus \widetilde{D} \rightarrow \mathbb{P}(d)\setminus D$.
We have $H_1(V^d \setminus \widetilde{D})\cong \mathbb{Z}$ (see \cite{D}Chapter 4, Corollary(1.4)).
Let $\gamma$ be the generator of $H_1(\mathbb{C}^*)$ represented by the loop
$\gamma(t)=e^{2\pi \sqrt{-1}t}, 0\le t\le 1$.
By Proposition \ref{prop:4-6} we see that the image of $\gamma$, which is represented by the loop
$t \mapsto e^{2\pi \sqrt{-1}t}\cdot F_0$, is $3(d-1)^2$ times a generator
of $H_1(V^d \setminus \widetilde{D})$.
On the other hand the loop
$$t \mapsto \left(\begin{array}{ccc}
e^{2\pi \sqrt{-1}t} & 0 & 0 \\
0 & e^{2\pi \sqrt{-1}t} & 0 \\
0 & 0 & e^{2\pi \sqrt{-1}t} 
\end{array}\right),\ 0\le t \le 1$$
in $GL(3;\mathbb{C})$, representing 3 times a generator of $H_1(GL(3;\mathbb{C}))\cong \mathbb{Z}$, is mapped 
to the loop $t \mapsto (e^{2\pi \sqrt{-1}t})^{-d} \cdot F_0$ by $\lambda$.
Hence in the commutative diagram
$$\xymatrix{
\mathbb{Z} \cong H_1(GL(3;\mathbb{C})) \ar[r]^{\lambda_*} \ar[d] & H_1(V^d \setminus \widetilde{D}) \ar[d] \\
\mathbb{Z}/3\mathbb{Z} \cong H_1(PGL(3)) \ar[r]^{\bar{\lambda}_*} & H_1(\mathbb{P}(d) \setminus D) \\
}$$
we have $\lambda_*(1)=\pm d(d-1)^2 \in \mathbb{Z} \cong H_1(V^d \setminus \widetilde{D})$  
so we can conclude $\bar{\lambda}_*(1\  {\rm mod}\ 3)=\pm d(d-1)^2\  {\rm mod}\ 3(d-1)^2$.
This completes the proof of Theorem \ref{thm:4-2}.

\section{The value of the Meyer function}
By Proposition \ref{prop:4-6}, we have $H^1(\pi_1(\mathbb{P}(d)\setminus D);\mathbb{Q})=0$.
Therefore $\bar{\phi}^d:=\phi^d \circ \chi_2$ is the unique 1-cochain
of $\pi_1(\mathbb{P}(d)\setminus D)$
satisfying $(\rho \circ \chi_2)^*\tau_g=\delta \bar{\phi}^d$.
In this section we will compute the value of $\bar{\phi}^d$ on a special element in
$\pi_1(\mathbb{P}(d)\setminus D)$ so called \textit{lasso}.

We first explain what a lasso is. Let $M$ be a connected complex manifold of dimension $m$ and
$N$ an irreducible hypersurface of $M$. Then the inclusion $M\setminus N \hookrightarrow M$ induces
the following exact sequence:
$$1 \rightarrow <\sigma> \rightarrow \pi_1(M\setminus N) \rightarrow \pi_1(M) \rightarrow 1.$$
Here $<\sigma>$ denotes the normal closure of an element $\sigma$ of $\pi_1(M\setminus N)$, which
is described in the following. Let $p$ be a non-singular point of $N$ and $(z_1,\ldots,z_m)$
a local coordinate system of $M$ around $p$ such that $N$ is defined by $z_1=0$.
For a sufficiently small $\varepsilon >0$, consider a loop defined in this coordinate system by
$$[0,1]\rightarrow M\setminus N,\ t\mapsto (\varepsilon e^{2\pi\sqrt{-1}t},0,\ldots,0)$$
based at $q=(\varepsilon,0,\ldots,0)$.
Joining this loop with a path from the base point of $M\setminus N$ to $q$, we get an
element $\sigma$ of $\pi_1(M\setminus N)$. Since $N$ is irreducible, the conjugacy class of
$\sigma$ in $\pi_1(M\setminus N)$ is independent of choices of $p$ and a local coordinate system.
Each element of this conjugacy class is called \textit{a lasso} around $N$.

Returning to $\pi_1(\mathbb{P}(d)\setminus D)$, $D$ is an irreducible hypersurface of $\mathbb{P}(d)$.
Let $\sigma^d \in \pi_1(\mathbb{P}(d)\setminus D)$ be a lasso around $D$. 
Since $\bar{\phi}^d$ is a class function (see Lemma \ref{lem:8-2} in Appendix), the values of $\bar{\phi}^d$ on
the conjugacy class of $\sigma^d$ is constant.

\begin{prop}For $d\ge 3$,
$$\bar{\phi}^d(\sigma^d)=-\frac{d+1}{3(d-1)}.$$
\label{prop:5-1}
\end{prop}
\begin{proof}
Choose $a^0$ and $Q \in U$ as in Lemma \ref{lem:4-5}.
In the proof of Proposition \ref{prop:4-6} we see that $Q$ meets $D$ transversely in $3(d-1)^2$ points.
Let $Q \cap D=\{ q_1,\ldots,q_{3(d-1)^2} \}$ and let $D_i$ ($i=1,\ldots,3(d-1)^2$) be a small closed
2-disk in $Q$ such that $q_i \in {\rm Int}D_i$ and $D_i \cap D_j=\emptyset$ for $i\neq j$.
We fix a base point of $Q_0:=Q\setminus \bigcup_{i=1}^{3(d-1)^2} {\rm Int}D_i$ and 
for each $i=1,\ldots,3(d-1)^2$, choose a based loop $\sigma_i$ in 
$Q_0$ such that $\sigma_i$ is free homotopic to the loop traveling once the boundary $\partial D_i$
by counter clockwise manner. Note that regarded as an element in $\pi_1(\mathbb{P}(d)\setminus D)$,
$\sigma_i$ is a lasso around $D$ hence we have $\bar{\phi}^d(\sigma_i)=\bar{\phi}^d(\sigma^d)$.

Let $g\colon M \rightarrow Q$ be as in the proof of Proposition \ref{prop:4-6} and set 
$M_0:=g^{-1}(Q_0)$ and $M_i:=g^{-1}(D_i), i=1,\ldots,3(d-1)^2$.
By Meyer's signature formula (\cite{Me}Satz 1) and the equation
$(\rho \circ \chi_2)^*\tau_g=\delta \bar{\phi}^d$, we obtain
$${\rm Sign}(M_0)=\sum_{i=1}^{3(d-1)^2}\bar{\phi}^d(\sigma_i)=3(d-1)^2\bar{\phi}^d(\sigma^d).$$
Since the topological type of $g^{-1}(q_i)$ is Lefschetz singular fiber of type I, we have ${\rm Sign}(M_i)=0$.
We compute by the Novikov additivity and Lemma \ref{lem:4-5}
$$1-d^2={\rm Sign}(M)={\rm Sign}(M_0)+\sum_{i=1}^{3(d-1)^2}{\rm Sign}(M_i)=3(d-1)^2\bar{\phi}^d(\sigma^d).$$
This completes the proof.
\end{proof}

In the rest of this section we consider the remaining case $d=2$. Since $V^2$ is the set
of quadratic forms each element of $V^2$ can be expressed by a $3\times 3$ symmetric matrix $S$.
In this view point $V^2 \setminus \widetilde{D}$ is the space of non-singular symmetric matrices and
the action of $GL(3;\mathbb{C})$ on $V^2 \setminus \widetilde{D}$ is given by
$$A\cdot S={\hskip -3pt}\ ^t{\hskip -1pt}A^{-1}\cdot S\cdot A^{-1},\ A\in GL(3;\mathbb{C}).$$
Since this action is transitive and the isotropy group of the unit matrix is the
complex orthogonal group $O_3(\mathbb{C})=\{ A\in GL(3;\mathbb{C})\ ;\ {\hskip -3pt}\ ^t{\hskip -1pt}A\cdot A=I \}$,
we have
$$V^2 \setminus \widetilde{D}\cong GL(3;\mathbb{C})/O_3(\mathbb{C}).$$
Also we have
$$\mathbb{P}(2)\setminus D\cong PGL(3)/SO_3(\mathbb{C}),$$
where $SO_3(\mathbb{C})=\{ A\in O_3(\mathbb{C})\ ;\ {\rm det}A=1 \}$ is regarded as a subgroup of
$PGL(3)$ by the injection $SO_3(\mathbb{C})\hookrightarrow PGL(3)$ induced by the projection
$GL(3;\mathbb{C})\rightarrow PGL(3)$.
Therefore, we obtain
\begin{eqnarray*}
(\mathbb{P}(2)\setminus D)_{PGL(3)}&=&EPGL(3)\times_{PGL(3)}(\mathbb{P}(2)\setminus D) \\
&\cong &EPGL(3)/SO_3(\mathbb{C})=BSO_3(\mathbb{C})\simeq BSO_3.
\end{eqnarray*}
The last homotopy equivalence holds because the natural inclusion $SO_3\hookrightarrow SO_3(\mathbb{C})$
is homotopy equivalence. In particular, we have
$$\Pi(2)\cong \pi_1(BSO_3)=1.$$

\section{The universal property of $(\mathbb{P}(d)\setminus D)_{PGL(3)}$}
In this section we will show the universal property of the space $(\mathbb{P}(d)\setminus D)_{PGL(3)}$.
In the latter part of the section, we consider the case $d=4$ more detail;
in particular, we prove that $\rho\colon \Pi(4)\rightarrow \Gamma_3$ is surjective.

We first make some definitions.
Let $\iota \colon X\rightarrow P$ be a continuous map and $h\colon P\rightarrow B$ a
$\mathbb{P}^2$ bundle whose structure group is $PGL(3)$.
We call $\xi=(X,\iota,P,h,B)$ \textit{a family of non-singular plane curves of degree $d$} if
\begin{enumerate}
\item $p:=h\circ \iota\colon X\rightarrow B$ is a continuous family of compact Riemann surfaces of genus
$g=\frac{1}{2}(d-1)(d-2)$, and
\item for each $b\in B$, the restriction $\iota|_{X_b}\colon X_b \rightarrow P_b$ is a holomorphic embedding
where $X_b=p^{-1}(b)$ and $P_b=h^{-1}(b)$.
\end{enumerate}
For each $b\in B$, the image $\iota(X_b)\subset P_b\cong \mathbb{P}^2$ is a non-singular plane curve of degree $d$.
Two such families $\xi_i=(X^i,\iota_i,P^i,h_i,B) , i=0,1$, are called \textit{isotopic} if there exists a family of
non-singular curves of degree $d$ over $B\times [0,1]$, denoted by
$\tilde{\xi}=(\widetilde{X},\tilde{\iota},\widetilde{P},\tilde{h},B\times [0,1])$, such that
for $i=0,1$, the restriction of $\tilde{\xi}$ to $B\times \{ i\}$ is isomorphic to $\xi_i$, i.e.,
for $i=0,1$, there exists a homeomorphism $\Psi_i\colon P^i \rightarrow \widetilde{P}|_{B\times\{i\}}$
and $\psi_i\colon X^i \rightarrow \widetilde{X}|_{B\times\{i\}}$ such that the diagram
$$\xymatrix{
X^i \ar[d]_{\iota_i} \ar[r]^{\psi_i} & \widetilde{X}|_{B\times\{i\}} \ar[d]^{\tilde{\iota}} \\
P^i \ar[d]_{h_i} \ar[r]^{\Psi_i} & \widetilde{P}|_{B\times\{i\}} \ar[d]^{\tilde{h}} \\
B \ar[r] & B\times \{i\},
}$$
where the last horizontal arrow is the homeomorphism $B\rightarrow B\times \{i\}$ given by $b\mapsto (b,i)$,
commutes and $\Psi_i$ (resp. $\psi_i$)  maps each fiber $P^i_b$ (resp. $X^i_b$) onto $\tilde{h}^{-1}(b,i)$
(resp. $(\tilde{h}\circ \tilde{\iota})^{-1}(b,i)$) biholomorphically.

For a given space $B$, we denote by $\mathcal{PC}_d(B)$ the set of all isotopy classes of
families of non-singular plane curves of degree $d$ over $B$. $\mathcal{PC}_d(\bullet)$ is contravariant;
for a given continuous map $f\colon B^{\prime}\rightarrow B$ we have a natural map
$\mathcal{PC}_d(B)\rightarrow \mathcal{PC}_d(B^{\prime})$ which assigns the isotopy class of
$\xi$ the isotopy class of the pull back of $\xi$ by $f$, which will be denoted by $f^*\xi$.
In fact, the isotopy class of $f^*\xi$ is uniquely determined by the homotopy class $[f]\in [B^{\prime},B]$.

Among such families of non-singular plane curves of degree $d$, there is a universal one.
Consider the inclusion map $\mathcal{F} \hookrightarrow (\mathbb{P}(d) \setminus D)\times \mathbb{P}^2$
and the first projection $(\mathbb{P}(d) \setminus D)\times \mathbb{P}^2 \rightarrow \mathbb{P}(d)\setminus D$.
For simplicity, we write $Y$ instead of $(\mathbb{P}(d) \setminus D)\times \mathbb{P}^2$.
Since these maps are $PGL(3)$-equivariant, we obtain
$$\iota_u\colon \mathcal{F}_{PGL(3)}\rightarrow Y_{PGL(3)} $$
and
$$h_u\colon Y_{PGL(3)}\rightarrow (\mathbb{P}(d)\setminus D)_{PGL(3)}.$$
The map $p_u:=h_u \circ \iota_u$ is the same as the map defined in section 1 and
$$\xi_u:=\left(\mathcal{F}_{PGL(3)},\iota_u,Y_{PGL(3)},h_u,(\mathbb{P}(d)\setminus D)_{PGL(3)}\right)$$
is a family of non-singular plane curves of degree $d$. The next theorem says that
$(\mathbb{P}(d)\setminus D)_{PGL(3)}$ is the classifying space for the functor $\mathcal{PC}_d(\bullet)$
and $\xi_u$ is the universal family.

\begin{thm}
\label{thm:6-1}
For any space $B$, the map
$$\eta \colon [B,(\mathbb{P}(d)\setminus D)_{PGL(3)}]\rightarrow \mathcal{PC}_d(B)$$
which assigns the homotopy class of $f\colon B\rightarrow (\mathbb{P}(d)\setminus D)_{PGL(3)}$
the isotopy class of the pull back $f^*\xi_u$, is bijective.
\end{thm}

In the following we shall construct the inverse of $\eta$.

Let $\xi=(X,\iota,P,h,B)$ be given. We divide the argument in three steps.

\noindent \textbf{Step 1.}\ We first consider the case when $P$ is trivial:
suppose that $P=B\times \mathbb{P}^2$. Then for each $b\in B$,
$\iota(X_b)\subset \{b\}\times \mathbb{P}^2\cong \mathbb{P}^2$ is a non-singular plane curve of
degree $d$, so the defining equation of $\iota(X_b)$ in $\mathbb{P}^2$ is uniquely determined as
an element of $\mathbb{P}(d)\setminus D$. Denoting it by $Eq(b)$, we obtain a map
$$Eq\colon B\rightarrow \mathbb{P}(d)\setminus D.$$

\begin{lem}The map $Eq$ is continuous.
\end{lem}
\begin{proof}
Regard $\mathbb{P}^2$ as the set of all complex lines
through the origin in $\mathbb{C}^3$. Then the holomorphic line bundle $\mathcal{O}(d)$ over $\mathbb{P}^2$
is given by
$$\mathcal{O}(d)=\mathcal{O}(1)^{\otimes d}=\bigcup_{\ell \in \mathbb{P}^2}{\rm Hom}(\ell,\mathbb{C})^{\otimes d}.$$
Let $p_2\colon B\times \mathbb{P}^2 \rightarrow \mathbb{P}^2$ be the second projection and
consider the pull back $L:=(p_2 \circ \iota)^*\mathcal{O}(d)$. $L\rightarrow X$ is a continuous
family over $B$ of holomorphic vector bundles. Now $H^0(\mathbb{P}^2;\mathcal{O}(d))$ is canonically
isomorphic to $V^d$ and for each $b\in B$ there is the natural homomorphism
$$\sigma_b\colon V^d\cong H^0(\mathbb{P}^2;\mathcal{O}(d)) \rightarrow
H^0\left(\iota(X_b);\mathcal{O}(d)|_{\iota(X_b)}\right)\cong H^0(X_b;L_b),$$
where $L_b$ is the restriction of $L$ to $X_b$.
Combining all $\sigma_b, b\in B$ together, we obtain a homomorphism of vector bundles
$$\sigma\colon B\times V^d \rightarrow \bigcup_{b\in B}H^0(X_b;L_b).$$
We see that for each $b\in B$, $\sigma_b$ is surjective and its kernel is 1-dimensional
generated by the defining equation of $\iota(X_b)$, i.e., $Eq(b)={\rm ker}\sigma_b$.
This shows that $Eq$ is continuous.
\end{proof}

We also define a $PGL(3)$-equivariant continuous map
$\Psi\colon B\times PGL(3)\rightarrow \mathbb{P}(d)\setminus D$ by
$$\Psi(b,g)=g\cdot Eq(b).$$
Here we regard $B\times PGL(3)$ as the trivial principal $PGL(3)$ bundle with
\textit{left} $PGL(3)$ action.

\noindent \textbf{Step 2.}\ We next consider the general case $\xi=(X,\iota,P,h,B)$.
Let $\{U_i\}_{i\in I}$ be an open covering of $B$ trivializing $h\colon P \rightarrow B$:
There is an isomorphism
$\varphi_i \colon h^{-1}(U_i)\rightarrow U_i\times \mathbb{P}^2$ for each $i$
and a system of transition functions
$g_{ij}\colon U_i \cap U_j \rightarrow PGL(3)$ for each $(i,j)$ satisfying $U_i \cap U_j\neq \emptyset$,
such that
$$(\varphi_i \circ \varphi_j^{-1})(b,p)=(b,g_{ij}(b)\cdot p),\ b\in U_i\cap U_j,\ p\in \mathbb{P}^2.$$
As in step 1, we have a continuous map
$Eq^i\colon U_i\rightarrow \mathbb{P}(d)\setminus D$ and a $PGL(3)$-equivariant map
$\Psi^i \colon U_i\times PGL(3)\rightarrow \mathbb{P}(d)\setminus D$ for each $i$.
Let $Q(\xi)$ be a principal $PGL(3)$ bundle over $B$ associated to $h\colon P\rightarrow B$:
namely $Q(\xi)$ is constructed from the disjoint union
$\coprod_{i\in I}U_i\times PGL(3)$ by identifying
$(b,g)\in U_i\times PGL(3)$ with $(b,g\cdot g_{ij}(b))\in U_j\times PGL(3)$ where $b\in U_i\cap U_j$.
We have $g_{ij}(b)\cdot Eq^j(b)=Eq^i(b)$ for $b\in U_i \cap U_j$
because $g\cdot C_{a}=C_{g\cdot a}$ for $g\in PGL(3),\ a\in \mathbb{P}(d)\setminus D$.
Therefore piecing all $\Psi^i,\ i\in I$ together, we obtain a $PGL(3)$ equivariant map
$\Psi\colon Q(\xi)\rightarrow \mathbb{P}(d)\setminus D$ and a continuous map
$$\Psi_{PGL(3)}\colon Q(\xi)_{PGL(3)}\rightarrow (\mathbb{P}(d)\setminus D)_{PGL(3)}.$$
Note that $Q(\xi)$ and $\Psi$ are determined up to isomorphism over $B$.

\noindent \textbf{Step 3.}\ The natural map
$$T\colon Q(\xi)_{PGL(3)}=EPGL(3)\times_{PGL(3)}Q(\xi)\rightarrow
PGL(3)\backslash Q(\xi)\cong B$$ is a homotopy equivalence because this is an $EPGL(3)$ bundle.
Taking a homotopy inverse map $\zeta \colon B\rightarrow Q(\xi)_{PGL(3)}$ of $T$, we set
$$\theta([\xi]):=[\Psi_{PGL(3)}\circ \zeta].$$ Here $[\xi]$ denotes
the element of $\mathcal{PC}_d(B)$ represented by $\xi$ and $[\Psi_{PGL(3)}\circ \zeta]$
denotes the element of $[B,(\mathbb{P}(d)\setminus D)_{PGL(3)}]$ represented by $\Psi_{PGL(3)}\circ \zeta$.
It is easy to see that
$$\theta\colon \mathcal{PC}_d(B) \rightarrow  [B,(\mathbb{P}(d)\setminus D)_{PGL(3)}]$$
is well defined.\\

Before starting the proof of Theorem \ref{thm:6-1}, we describe the above construction
applied to the family $\xi_u$. In the below, $e\in EPGL(3),\ g,h\in PGL(3)$ and
$a\in \mathbb{P}(d)\setminus D$. We can write
\begin{equation}
\label{eq:4}
Q(\xi_u)\cong \left((\mathbb{P}(d)\setminus D)\times PGL(3) \right)_{PGL(3)}
\end{equation}
where the action of $PGL(3)$ on $(\mathbb{P}(d)\setminus D)\times PGL(3)$ is diagonal, i.e.,
$$g\cdot (a,h)=(g\cdot a,g\cdot h),$$
and the left action of $PGL(3)$ on the right hand side of (\ref{eq:4}) is given by
$$g\cdot[e,(a,h)]=[e,(a,h\cdot g^{-1})].$$
The $PGL(3)$-equivariant map
$\Psi_u\colon Q(\xi_u)\rightarrow \mathbb{P}(d)\setminus D$ defined as in step 2 is given by
$$\Psi_u([e,(a,g)])=g^{-1}\cdot a,$$
and moreover, the induced map
$Q(\xi_u)_{PGL(3)}\rightarrow (\mathbb{P}(d)\setminus D)_{PGL(3)}$ has a section $s_u$
given by
$$s_u([e,a])=[e,[e,(a,1)]].$$

\noindent{\it Proof of Theorem \ref{thm:6-1}.}
We first prove $\eta \circ \theta=id_{\mathcal{PC}_d(B)}$. Let $\xi=(X,\iota,P,h,B)$ be given.
By construction, there is the canonical isomorphism $T^*\xi \rightarrow \Psi_{PGL(3)}^*\xi_u$
as families of non-singular plane curves of degree $d$ over $Q(\xi)_{PGL(3)}$. Thus we have
$$(\Psi_{PGL(3)}\circ \zeta)^*\xi_u=\zeta^*\Psi_{PGL(3)}^*\xi_u
=\zeta^*T^*\xi=(T\circ \zeta)^*\xi.$$
Since $T\circ \zeta$ is homotopic to the identity map of $E$, this shows
$\eta \circ \theta=id_{\mathcal{PC}_d(B)}$.

We next show $\theta \circ \eta=id_{[B,(\mathbb{P}(d)\setminus D)_{PGL(3)}]}$. Let
$f\colon B\rightarrow (\mathbb{P}(d)\setminus D)_{PGL(3)}$ be a continuous map.

Starting from the family $f^*\xi$ and tracing the construction of
$\theta$, we construct the map
$$\Psi_{PGL(3)} \colon Q(f^*\xi_u)_{PGL(3)}\rightarrow (\mathbb{P}(d)\setminus D)_{PGL(3)}.$$
$Q(f^*\xi_u)_{PGL(3)}$ is naturally isomorphic to the pull back of
$Q(\xi_u)_{PGL(3)}\rightarrow (\mathbb{P}(d)\setminus D)_{PGL(3)}$ by $f$.
Thus pulling back the section $s_u$, we obtain a map
$\zeta^{\prime}:=f^*s_u\colon B\rightarrow Q(f^*\xi_u)_{PGL(3)}$
such that $T \circ \zeta^{\prime}=id_B$ and $\Psi_{PGL(3)} \circ \zeta^{\prime}=f$.
Then $\zeta^{\prime}$ is a homotopy inverse of $T$ and
$\theta \circ \eta([f])=[\Psi_{PGL(3)} \circ \zeta^{\prime}]=[f]$, so we obtain
$\theta \circ \eta=id_{[B,(\mathbb{P}(d)\setminus D)_{PGL(3)}]}$. \qed

We call any representative of $\theta([\xi])$
\textit{the classifying map for the family $\xi$}.

For $d=4$, we don't have to consider $\mathbb{P}^2$ bundles.
Recall that a non-hyperelliptic Riemann surface $C$ of genus 3 can be realized as
a non-singular plane curve of degree 4 (=plane quartic) by the canonical embedding.
This means that the canonical map
$$\iota_C\colon C\rightarrow \mathbb{P}(H^0(C;K_C)^{\vee})$$
where $H^0(C;K_C)$ is the space of holomorphic 1-forms on $C$, is an embedding
and if we identify $\mathbb{P}(H^0(C;K_C)^{\vee})$ with $\mathbb{P}^2$ by a choice
of a basis of $H^0(C;K_C)$, the image of $C$ is a non-singular plane curve of degree 4.
The defining equation of the image is uniquely determined as an element of $\mathbb{P}(4)\setminus D$.

Let $p\colon X\rightarrow B$ be a continuous family of compact Riemann surfaces of genus 3.
We call it \textit{a non-hyperelliptic family of genus 3} if the complex structure of each fiber
$p^{-1}(b),b\in B$ is non-hyperelliptic. Two such families $p_i\colon X_i \rightarrow B \ (i=0,1)$
are called isotopic if there exists a non-hyperelliptic family of genus 3 
over $B\times [0,1]$ such that for $i=0,1$, its restriction to $B\times \{i\}$ is isomorphic to 
$p_i\colon X_i\rightarrow B$ as continuous family of Riemann surfaces over $B\cong B\times \{i\}$.

For a given space $B$, we denote by $\mathcal{NH}_3(B)$ the set of all isotopy classes of
non-hyperelliptic families of genus 3 over $B$. Then the forgetful functor
\begin{equation}
\label{eq:5}
\mathcal{PC}_4(\bullet) \rightarrow \mathcal{NH}_3(\bullet)
\end{equation}
defined by an obvious manner, is bijective. For, let $p\colon X\rightarrow B$ be a given
non-hyperelliptic family of genus 3. Set $\Lambda_X:=\bigcup_{b\in B}H^0(p^{-1}(b);K_b)$,
where $H^0(p^{-1}(b);K_b)$ denotes the space of holomorphic 1-forms on $p^{-1}(b)$.
This has the structure of complex vector bundle over $B$. Projectivising the dual of $\Lambda_X$,
we obtain a $\mathbb{P}^2$ bundle
$$h^{\prime}\colon P^{\prime}=\bigcup_{b\in B}\mathbb{P}(H^0(p^{-1}(b);K_b)^{\vee})\rightarrow B,$$
and piecing the fiberwise canonical maps $\iota_{X_b},\ b\in B$ together, we get
a map $\iota\colon X\rightarrow P^{\prime}$. Then we obtain an element
$\xi=(X,\iota,P^{\prime},h^{\prime},B)\in {PC}_4(\bullet)$. 
This correspondence gives the inverse of (\ref{eq:5}). 

We continue the consideration of the case $d=4$. We next prove that:
\begin{prop}
\label{prop:6-3}
The homomorphism $\rho \colon \Pi(4)\rightarrow \Gamma_3$ is surjective.
\end{prop}

Combining this with Theorem \ref{thm:4-1}, which implies that $\rho^*\colon H^2(\Gamma_3;\mathbb{Q})\rightarrow
H^2(\Pi(4);\mathbb{Q})$ is not injective, we see that the order of the kernel of $\rho$ is infinite.

\noindent{\it Proof of Proposition \ref{prop:6-3}.}
Let $\mathcal{T}_3$ be the Teichm\"uller space of compact Riemann surfaces of genus 3
and $H_3$ the hyperelliptic locus of $\mathcal{T}_3$; namely the set of marked Riemann surfaces
whose complex structure is hyperelliptic. $H_3$ is a complex analytic closed submanifold
of codimension 1 with infinitely many components (see \cite{N}p.259-260).
In particular, $\mathcal{T}_3 \setminus H_3$ is path connected. 

We recall; there is a holomorphic family $\pi\colon V_3\rightarrow \mathcal{T}_3$
called the universal Teichm\"uller curve, whose fiber over the marked Riemann surface
$[f,C]$ is isomorphic to $C$;
the mapping class group $\Gamma_3$ acts on $\mathcal{T}_3$ and $V_3$, and
$\pi$ is equivariant with respect to these actions;
it is well known that the quotient space $\Gamma_3 \backslash \mathcal{T}_3$
is the Riemann moduli space.
Since the action of $\Gamma_3$ on $\mathcal{T}_3$ preserves $H_3$,
$\Gamma_3$ also acts on $\mathcal{T}_3 \setminus H_3$ and its inverse
image by $\pi$. Restricting $\pi$ to
$\mathcal{T}_3 \setminus H_3$ and taking the Borel construction, we obtain a
non-hyperelliptic family of genus 3 over $(\mathcal{T}_3 \setminus H_3)_{\Gamma_3}$.

It is not difficult to see that this family also have the universal property
which the family $p_u$ over $(\mathbb{P}(4)\setminus D)_{PGL(3)}$ has.
Therefore, $(\mathcal{T}_3 \setminus H_3)_{\Gamma_3}$ is homotopy equivalent
to $(\mathbb{P}(4)\setminus D)_{PGL(3)}$ hence its fundamental group is isomorphic
to $\Pi(4)$.

By the homotopy exact sequence of the $\mathcal{T}_3\setminus H_3$ bundle
$(\mathcal{T}_3 \setminus H_3)_{\Gamma_3} \rightarrow B\Gamma_3=K(\Gamma_3,1)$ we
obtain an exact sequence
$$\Pi(4)\cong \pi_1((\mathcal{T}_3 \setminus H_3)_{\Gamma_3})
\stackrel{\rho^{\prime}}{\rightarrow} \pi_1(B\Gamma_3)=\Gamma_3
\rightarrow \pi_0(\mathcal{T}_3 \setminus H_3).$$

We notice that the homomorphism $\rho^{\prime}$ just coincides with the topological monodromy over
$(\mathcal{T}_3 \setminus H_3)_{\Gamma_3}$, and $\pi_0(\mathcal{T}_3 \setminus H_3)$ is one point.
This shows $\rho^{\prime}$ is surjective, so $\rho$ is. \qed

\section{Local signature for 4-dimensional\\
non-hyperelliptic fibration of genus 3}
As an application, we will define the local signature for the set of all fiber germs of
4-dimensional fiber spaces whose general fibers are non-hyperelliptic Riemann surfaces
of genus 3, using the Meyer function $\phi^4$.
This local signature is used to derive a signature formula for a class of 4-dimensional fiber spaces,
whose general fibers are non-hyperelliptic Riemann surfaces of genus 3.

Let $\Delta$ be a closed oriented 2-disk and $p$ its center. A 4-tuple
$\mathcal{F}=(E,\pi,\Delta,p)$ is called \textit{a fiber germ of non-hyperelliptic family of genus 3} if
\begin{enumerate}
\item $E$ is a $C^{\infty}$ manifold of dimension 4 and $\pi \colon E\rightarrow \Delta$ is a $C^{\infty}$ map,
\item the restriction of $\pi$ to $\Delta \setminus \{p \}$ is a non-hyperelliptic family of genus 3.
\end{enumerate}
Note that $E$ has the natural orientation and compact, hence its signature ${\rm Sign}(E)$ is defined.
Two such germs $(E,\pi,\Delta,p)$ and $(E^{\prime},\pi^{\prime},\Delta^{\prime},p^{\prime})$ are called
\textit{equivalent} if there exist a smaller disk $\Delta_0 \subset \Delta$
(resp. $\Delta^{\prime}_0 \subset \Delta^{\prime}$)
whose center is $p$ (resp. $p^{\prime}$), and there exist orientation preserving diffeomorphisms
$\varphi \colon (\Delta_0,p)\rightarrow (\Delta^{\prime}_0,p^{\prime})$
and $\tilde{\varphi}\colon \pi^{-1}(\Delta_0)\rightarrow \pi^{\prime -1}(\Delta^{\prime}_0)$ such that
$\varphi \circ \pi=\pi^{\prime} \circ \tilde{\varphi}$ and
$$\tilde{\varphi}|_{\pi^{-1}(\Delta_0 \setminus \{ p\})}\colon
\pi^{-1}(\Delta_0 \setminus \{ p\})\rightarrow
\pi^{\prime -1}(\Delta^{\prime}_0 \setminus \{ p^{\prime}\})$$ maps each fiber biholomorphically.

Let $\mathcal{NH}_3$ denote the set of all equivalence classes of such 4-tuples.
We denote the element of $\mathcal{NH}_3$ also by $\mathcal{F}=(E,\pi,\Delta,p)$.
For $\mathcal{F}=(E,\pi,\Delta,p)\in \mathcal{NH}_3$, $\gamma$ denotes the element of
$\pi_1(\Delta \setminus \{ p\} )$ traveling once the boundary $\partial \Delta$ by counter clockwise manner.
We denote by $\mathcal{F}^0$ the restriction of $\pi\colon E \rightarrow \Delta$ to $\Delta\setminus \{ p\}$.
$\mathcal{F}^0$ is a non-hyperelliptic family of genus 3 and can be considered as an element of
$\mathcal{PC}_4(\Delta\setminus \{ p\})$ in view of (\ref{eq:5}).

\begin{dfn}Define ${\rm loc.sig}^{\mathcal{Q}}\colon \mathcal{NH}_3 \rightarrow \mathbb{Q}$ by
$${\rm loc.sig}^{\mathcal{Q}}(\mathcal{F}):=\phi^4(\theta(\mathcal{F}^0)_*(\gamma))+{\rm Sign}(E).$$
\end{dfn}
Here, $\theta(\mathcal{F}^0)_*$ is the homomorphism from $\pi_1(\Delta \setminus \{p \})$ to
$\Pi(4)$ induced by the classifying map $\theta(\mathcal{F}^0)$ for $\mathcal{F}^0$.
It is assumed that suitable base points of $\Delta \setminus \{ p\}$ and
$(\mathbb{P}(4)\setminus D)_{PGL(3)}$ are chosen. Since $\phi^4$ is a class function,
we don't have to care about base point so we omit it.

We call a triple $(E,\pi,B)$ \textit{a 4-dimensional non-hyperelliptic fibration of genus 3} if
\begin{enumerate}
\item $E$ (resp. $B$) is a closed oriented $C^{\infty}$-manifold of dimension 4 (resp. 2) and
$\pi\colon E\rightarrow B$ is a $C^{\infty}$-map,
\item there exist finitely many points $b_1,\ldots,b_n \in B$ such that the restriction of
$\pi$ to $B\setminus \{ b_1,\ldots,b_n\} $ is a non-hyperelliptic family of genus 3. 
\end{enumerate}
For $i=1,\ldots,n$, we obtain an element of $\mathcal{NH}_3$ by restricting $\pi$ to a small
closed disk neighborhood of $b_i$. we denote it by $\mathcal{F}_i$. Then, we obtain
\begin{thm}[The signature formula]
Let $(E,\pi,B)$ be a 4-dimensional non-hyperelliptic fibration of genus 3. Then
$${\rm Sign}(E)=\sum_{i=1}^n{\rm loc.sig}^{\mathcal{Q}}(\mathcal{F}_i).$$
\label{thm:7-2}
\end{thm}
\begin{proof}
For $i=1,\ldots,n$, take a small closed 2-disk $D_i$ around $b_i$
so that they don't intersect each other. Then $\mathcal{F}_i=(\pi^{-1}(D_i),\pi,D_i,b_i)$.
We denote by $\mathcal{F}^0_i$ the restriction of $\pi$ to $D_i \setminus \{ b_i\}$ and set
$B_0:=B \setminus \bigcup^{n}_{i=1}{\rm Int}D_i$. By Meyer's signature formula, we get
$${\rm Sign}(\pi^{-1}(B_0))=\sum^{n}_{i=1}\phi^4(\theta(\mathcal{F}^0_i)_*(\gamma)).$$
Using the Novikov additivity, we compute
\begin{eqnarray*}
{\rm Sign}(E)&=& {\rm Sign}(\pi^{-1}(B_0))+\sum^{n}_{i=1}{\rm Sign}(\pi^{-1}(D_i)) \\
&=& \sum^{n}_{i=1}\phi^4(\theta(\mathcal{F}^0_i)_*(\gamma))+\sum^{n}_{i=1}{\rm Sign}(\pi^{-1}(D_i))\\
&=& \sum^{n}_{i=1}{\rm loc.sig}^{\mathcal{Q}}(\mathcal{F}_i).
\end{eqnarray*}
\end{proof}

\begin{cor}
Let $g\colon E\rightarrow B$ be
a non-hyperelliptic family of genus 3 over a closed oriented surface $B$.
Then ${\rm Sign}(E)=0.$
\end{cor}

We compute some examples. Comparing the following computations with those in
\cite{AK} and \cite{Y}, we see that their values coincide.

\noindent \textbf{Singular fiber of type I.}\ Let $\Delta \subset \mathbb{P}(4)$ be a
closed 2-disk intersecting with $D$ only in its center $p\in \Delta$ transversely.
Let $\pi_I\colon E_I\rightarrow \Delta$ be the restriction of $\bar{\mathcal{F}}\rightarrow \mathbb{P}(4)$
to $\Delta$. Then $E_I$ is smooth by Lemma \ref{lem:2-4} and $\mathcal{F}_I=(E_I,\pi_I,\Delta,p)$ is
a fiber germ of non-hyperelliptic family of genus 3. By Lemma \ref{lem:2-3} the
topological type of $\pi_I^{-1}(p)$ is Lefschetz singular fiber of type I, therefore
we also call $\mathcal{F}_I\in \mathcal{NH}_3$ \textit{a singular fiber germ of type I}.
The signature of $E_I$ is $0$ and by definition, the inclusion
$\Delta \setminus \{ p\} \hookrightarrow\mathbb{P}(4)\setminus D \hookrightarrow (\mathbb{P}(4)\setminus D)_{PGL(3)}$
is the classifying map for $\mathcal{F}_I^0$ and the boundary of $\Delta$ is a lasso about $D$.
Therefore, by Proposition \ref{prop:5-1}, we have

\begin{prop}
$${\rm loc.sig}^{\mathcal{Q}}(\mathcal{F}_I)=-\frac{5}{9}.$$
\label{prop:7-4}
\end{prop}

\ \\
\noindent \textbf{Hyperelliptic fiber.}\ Let $F\in V^4\setminus \{0\}$ be a polynomial such that $C_F$
intersects with the non-singular conic $C\colon yz-x^2=0$ in $8$ points, and
let $\Delta$ be a small closed 2-disk around $0\in \mathbb{C}$ with the complex coordinate $s$.
Let $S_F$ be the hypersurface in $\Delta \times \mathbb{P}^2$ defined by the equation
$$(yz-x^2)^2+s^2F(x,y,z)=0.$$
$S_F$ is singular along $C^{\prime}=\{ 0\} \times C$.
Blowing up $\Delta \times \mathbb{P}^2$ along $C$, let $\widetilde{S_F}$ be the
proper transform of $S_F$ and $\pi\colon \widetilde{S_F}\rightarrow \Delta$ the composition
of $\widetilde{S_F}\rightarrow S_F$ and the first projection $S_F\rightarrow \Delta$.
Then $\widetilde{S_F}$ is non-singular and the exceptional divisor $\pi^{-1}(0)$
is a non-singular hyperelliptic curve of genus 3 with a natural projection
onto $C^{\prime}\cong \mathbb{P}^1$, which is a double cover.

Choose $\Delta$ small enough so that the singular fiber
of $\pi$ is $\pi^{-1}(0)$ only.
Set $\mathcal{F}_h=(\widetilde{S_F},\pi,\Delta,0)$ and call this fiber germ
\textit{a hyperelliptic germ}.
Let $\ell_h$ be the corresponding loop in $\mathbb{P}(4)\setminus D$ defined by
$$\ell_h(t)=(yz-x^2)^2+(\varepsilon e^{2\pi \sqrt{-1}t})^2F(x,y,z),\ 0\le t\le 1,$$
where $\varepsilon$ is the radius of $\Delta$.

\begin{prop}
$${\rm loc.sig}^{\mathcal{Q}}(\mathcal{F}_h)=\bar{\phi}^4([\ell_h])=\frac{4}{9}.$$
\end{prop}
\begin{proof}
We first note that
${\rm loc.sig}^{\mathcal{Q}}(\mathcal{F}_h)=\bar{\phi}^4([\ell_h])$ since 
a hyperelliptic germ is topologically trivial.

The set $W$ of all polynomials in $V^4$ such
that the corresponding curve intersects with $C$ in $8$ points is a non-empty Zariski open subset of $V^4$.
Since $[\ell_h]$ and ${\rm Sign}(\widetilde{S_F})$ does not change under any small
perturbation of $F$ in $V^4$, it suffices to show the proposition for a particular element
of $W$. But by the same reason as in Lemma \ref{lem:4-4}, there is actually an element $F\in W$
such that the map
$$\mathbb{P}^1 \rightarrow \mathbb{P}(4),\ [w_0\colon w_1]\mapsto w_0^2(yz-x^2)^2+w_1^2F(x,y,z),$$
does not meet $E$ and is transverse to $D$, except at $[w_0\colon w_1]=[1\colon 0]$.
Then for this choice of $F$, the complex surface $S$ in $\mathbb{P}^1\times \mathbb{P}^2$
defined by the equation
$$w_0^2(yz-x^2)^2+w_1^2F(x,y,z)=0,$$
has singularities only along the conic $\{[1\colon 0]\} \times C$. After blowing up
$\mathbb{P}^1\times \mathbb{P}^2$ along this conic, we obtain the proper transform
$\widetilde{S}$ of $S$. By the choice of $F$, $\widetilde{S}$ is non-singular.
The composition of $\widetilde{S}\rightarrow S$ and the first projection $S\rightarrow \mathbb{P}^1$
is a family of algebraic curves whose all singular fiber germs are singular fiber germ of type I
except the fiber germ around $[1\colon 0]$, and the fiber germ around $[1\colon 0]$ is a hyperelliptic germ.
The invariants of $\widetilde{S}$ are computed as: $c_1^2(\widetilde{S})=-6,\ c_2(\widetilde{S})=18$,
and ${\rm Sign}(\widetilde{S})=-14$.

Now the number of singular fiber germs of type I is equal to the total Euler contribution
$$18-2(2-2\cdot 3)=26.$$
Note that a hyperelliptic germ, which is topologically trivial, does not contribute to the Euler number.
By Theorem \ref{thm:7-2} and Proposition \ref{prop:7-4}, we have
$$-14=-\frac{5}{9}\times 26+{\rm loc.sig}^{\mathcal{Q}}(\mathcal{F}_h),$$
hence ${\rm loc.sig}^{\mathcal{Q}}(\mathcal{F}_h)=\frac{4}{9}$.
\end{proof}

\ \\
\noindent \textbf{Singular fiber of type II.}\ Let $\Delta$ be as in the previous example, and
let $S$ be the surface in $\Delta \times \mathbb{P}^2$ defined by
$$z^3x+y^2x^2+y^4+s^6x^4=0.$$
$S$ has an isolated singularity at $p_0=(0,[1\colon 0\colon 0])$ so called a singularity of type $\widetilde{E}_8$.
The inverse image $C_2$ of $0\in \Delta$ by the
first projection $p_1 \colon S \rightarrow \Delta$ is a curve of geometric genus 2 with
one cusp singularity.

Let $\varpi:\widetilde{S}\rightarrow S$ be the minimal resolution of the singularity of $S$ at $p_0$.
Then the exceptional curve is a non-singular elliptic curve $C_1$ with self intersection number $-1$.
If $\Delta$ is small enough, $\mathcal{F}_{II}=(\widetilde{S},p_1 \circ \varpi,\Delta,0)$
is a fiber germ of non-hyperelliptic family of genus 3. The topological type of the singular fiber
$(p_1 \circ \varpi)^{-1}(0)$ is obtained by the disjoint union of $C_1$ and $C_2$ by identifying
a point of $C_1$ with the cusp singularity of $C_2$, that is, Lefschetz singular fiber of type II.
We call $\mathcal{F}_{II}$ \textit{a singular fiber germ of type II}.

Let $\ell_{II}$ be the corresponding loop in $\mathbb{P}(4)\setminus D$ defined by
$$\ell_{II}(t)=z^3x+y^2x^2+y^4+(\varepsilon e^{2\pi \sqrt{-1}t})^6x^4,\ 0\le t\le 1.$$

\begin{prop}
$${\rm loc.sig}^{\mathcal{Q}}(\mathcal{F}_{II})=\frac{1}{3},\hspace{1em}
\bar{\phi}^4([\ell_{II}])=\frac{4}{3}.$$
\end{prop}
\begin{proof}
In this case $\bar{\phi}^4([\ell_{II}])={\rm loc.sig}^{\mathcal{Q}}(\mathcal{F}_{II})+1$
because the intersection form of $\widetilde{S}$ is given by $\left( \begin{array}{cc}
0 & 1 \\ 1 & 0 \\ \end{array} \right)$ hence ${\rm Sign}(\widetilde{S})=-1$.

We perturb $S$ slightly by adding a higher term about $s$; consider
the surface in $\Delta \times \mathbb{P}^2$ defined by
$$z^3x+y^2x^2+y^4+s^6x^4+s^mF(x,y,z)=0,$$
where $m$ is an integer $\ge 7$ and $F$ is a polynomial in $V^4$.
The singularity of this surface remains at the origin and is still of type $\widetilde{E}_8$. Taking the minimal
resolution of this singularity and taking $\Delta$ to be smaller if needed,
we obtain a new fiber germ $\mathcal{F}_{II}^{\prime}$ and a new loop $\ell_{II}^{\prime}$ in 
$\mathbb{P}(4)\setminus D$. This perturbation does not influence the value of $\bar{\phi}^4$
and the topology of the fiber neighborhood of the singular fiber. So it suffices to compute
${\rm loc.sig}^{\mathcal{Q}}(\mathcal{F}_{II}^{\prime})$.

Let $S^{\prime}$ be the complex surface in $\mathbb{P}^1\times \mathbb{P}^2$
defined by the equation
$$w_0^m(z^3x+y^2x^2+y^4)+w_0^{m-6}w_1^6x^4+w_1^mF(x,y,z)=0,$$
and let $\widetilde{S^{\prime}}\rightarrow S^{\prime}$ be the minimal resolution
of the singularity of $S^{\prime}$ at $p_0=([1\colon 0],[1\colon 0\colon 0])$.
If a generic $F$ is chosen, then $\widetilde{S^{\prime}}$ is non-singular and
the singular fiber germs of the family of algebraic curves
$\widetilde{S^{\prime}}\rightarrow S^{\prime}\rightarrow \mathbb{P}^1$ are all of type I
except the fiber germ around $[1\colon 0]$, and the fiber germ around $[1\colon 0]$ is 
$\mathcal{F}_{II}^{\prime}$.
The invariants of $\widetilde{S^{\prime}}$ are computed as:
$c_1^2(\widetilde{S^{\prime}})=9m-17,\ c_2(\widetilde{S^{\prime}})=27m-19$,
and ${\rm Sign}(\widetilde{S^{\prime}})=-15m+7$.

Now the number of singular fiber germs of type I is equal to
$$27m-19-2(2-2\cdot 3)-1=27m-12.$$
This time, the fiber over $[1\colon 0]$ do contribute to the Euler number.

By the signature formula,
$$-15m+7=-\frac{5}{9}\times(27m-12)+{\rm loc.sig}^{\mathcal{Q}}(\mathcal{F}_{II}^{\prime}).$$
Thus we obtain ${\rm loc.sig}^{\mathcal{Q}}(\mathcal{F}_{II})=
{\rm loc.sig}^{\mathcal{Q}}(\mathcal{F}_{II}^{\prime})=\frac{1}{3}$.
\end{proof}

\section{Appendix}
In this appendix we give a definition of Meyer's signature cocycle in the
form used in the present paper and review its properties. For details, see
W.\ Meyer's original paper \cite{Me}.

We first explain the topological monodromy of surface bundles.
Let $\pi\colon E\rightarrow B$ be an oriented
$\Sigma_g$ bundle whose structure group is the group of all orientation
preserving diffeomorphisms of $\Sigma_g$. Choose a base point $b_0\in B$ and
fix an identification $\phi\colon \Sigma_g \stackrel{\cong}{\rightarrow}\pi^{-1}(b_0)$.
For each based loop $\ell\colon [0,1]\rightarrow B$ the pull back $\ell^*(E)\rightarrow [0,1]$
of $\pi\colon E\rightarrow B$ by $\ell$ is trivial. Hence there exist a trivialization
$\Phi\colon \Sigma_g \times [0,1]\rightarrow \ell^*(E)$ such that $\Phi(x,0)=\phi(x)$.
By assigning the isotopy class of $\Phi(x,1)^{-1}\circ \phi$ to the homotopy class of $\ell$,
we obtain a map $\chi\colon \pi_1(B,b_0)=\pi_1(B)\rightarrow \Gamma_g$.
This map becomes a homomorphism under the conventions;
1) for any two mapping classes $f_1$ and $f_2$, the multiplication $f_1\circ f_2$ means
that $f_2$ is applied first,
2) for any two homotopy classes of based loops $\ell_1$ and $\ell_2$, their product
$\ell_1\cdot \ell_2$ means that $\ell_1$ is traversed first.
$\chi$ is called \textit{the topological monodromy of $\pi\colon E\rightarrow B$} and determined
up to inner automorphisms of $\Gamma_g$.

Let $P$ denote the pair of pants, i.e., $P=S^2\setminus \bigcup_{i=1}^{3}{\rm Int}D_i$
where $D_i,i=1,2$, and 3 are the three disjoint closed disks in the 2-sphere $S^2$.
Choose a base point $p_0\in {\rm Int}P$ and fix a based loop $\ell_1$ and $\ell_2$
such that $\ell_i$ is free homotopic to the loop traveling once the boundary $\partial D_i$
by counter clockwise manner ($i=1,2$).
For $(f_1,f_2)\in \Gamma_g\times \Gamma_g$, we can construct an oriented $\Sigma_g$ bundle
$E(f_1,f_2)$ over $P$ such that the topological monodromy
$\chi \colon\pi_1(P)\rightarrow \Gamma_g$ sends $[\ell_i]$ to $f_i$ for $i=1,2$.
(If $g\ge 2$, the isomorphism class of this bundle is unique.)
$E(f_1,f_2)$ is a compact $C^{\infty}$-manifold of dimension 4 and has the natural orientation
induced by the orientation of $P$ and that of the fibers. Then the signature of $E(f_1,f_2)$
is defined and we set
$$\tau_g(f_1,f_2):=-{\rm Sign}(E(f_1,f_2)).$$
This turns out to be well defined even when $g=1$, and $\tau_g\colon \Gamma_g\times \Gamma_g \rightarrow \mathbb{Z}$
is called \textit{Meyer's signature cocycle}. The basic properties of $\tau_g$ are \\
$(1)\ \tau_g(f_1f_2,f_3)+\tau_g(f_1,f_2)=\tau_g(f_1,f_2f_3)+\tau_g(f_2,f_3)$; \\
$(2)\ \tau_g(f_1,1)=\tau_g(1,f_1)=\tau_g(f_1,f_1^{-1})=0$; \\
$(3)\ \tau_g(f_1^{-1},f_2^{-1})=-\tau_g(f_1,f_2)$; \\
$(4)\ \tau_g(f_1,f_2)=\tau_g(f_2,f_1)$; \\
$(5)\ \tau_g(f_3f_1f_3^{-1},f_3f_2f_3^{-1})=\tau_g(f_1,f_2)$, \\
where $f_1,f_2$, and $f_3$ are elements of $\Gamma_g$.

For an oriented $\Sigma_g$ bundle $\pi\colon E\rightarrow B$ and a choice of base point $b_0$ of $B$,
we obtain a 2-cocycle $\chi^*\tau_g$ of $\pi_1(B)=\pi_1(B,b_0)$ by pulling back $\tau_g$ by
the topological monodromy $\chi\colon \pi_1(B)\rightarrow \Gamma_g$. Although $\chi$ is
determined only up to conjugacy, $\chi^*\tau_g$ is uniquely determined by the property $(5)$ of $\tau_g$ above.
Moreover, $\chi^*\tau_g$ does not depend on the choice of base point of $B$ in the following
sense: suppose $b_0^{\prime}\in B$ and $b_0$ are in the same path component of $B$ then under
any isomorphism $\pi_1(B,b_0)\cong \pi_1(B,b_0^{\prime})$ using a path from $b_0$ to $b_0^{\prime}$,
two cocycles of $\pi_1(B,b_0)$ and $\pi_1(B,b_0^{\prime})$ defined as the pull back of $\tau_g$
by topological monodromies, correspond to each other.

Let $G$ be a group and $\varphi:G\rightarrow \Gamma_g$ a homomorphism.

\begin{dfn}
A $\mathbb{Q}$-valued 1-cochain $\phi \colon G \rightarrow \mathbb{Q}$ is called a Meyer function
with respect to the pull back $\varphi^*\tau_g$ of $\tau_g$ by $\varphi$
if it satisfies $\delta \phi =\varphi^*\tau_g$, i.e., $\phi$ cobounds the 2-cocycle $\varphi^*\tau_g$.
\end{dfn}

If a Meyer function exists on $G$, the cohomology class
$\varphi^*[\tau_g]\in H^2(G;\mathbb{Z})$ is torsion.
The following properties of $\phi$ are easily derived by the above properties of $\tau_g$
(see also \cite[Proposition 3.1]{E}).

\begin{lem}
\label{lem:8-2}
If $\phi$ is a Meyer function with respect to $\varphi^*\tau_g$, we have \\
$(1) \ \phi(xy)=\phi(x)+\phi(y)-\varphi^*\tau_g(x,y);$ \\
$(2) \ \phi(1)=0;$ \\
$(3) \ \phi(x^{-1})=-\phi(x);$ \\
$(4) \ \phi(yxy^{-1})=\phi(x),$ \\
where $x,y \in G.$
\end{lem}

\noindent \textbf{Acknowledgements.}\ I would like to express special thanks to my
advisor, Nariya Kawazumi, for many comments, proofreading, and offering me the opportunity to study
the work of W.Meyer and localization of the signature.
My thanks also go to Tadashi Ashikaga, for helpful suggestions on computations of the local signature
performed in section $7$ of this paper. This research is supported by JSPS Research
Fellowships for Young Scientists (19$\cdot$5472).

\noindent \textsc{Yusuke Kuno\\
Graduate School of Mathematical Sciences,\\
The University of Tokyo,\\
3-8-1 Komaba Meguro-ku Tokyo 153-0041, JAPAN}

\noindent \texttt{E-mail address:kunotti@ms.u-tokyo.ac.jp}

\end{document}